\documentclass[12pt,twoside]{amsart}
\usepackage{euler}
\usepackage{epsfig}
\usepackage{latexsym}
\usepackage{graphicx}
\usepackage{enumerate}

 \def\Dj{\hbox{D\kern-.73em\raise.30ex\hbox{-} \raise-.30ex\hbox{}}}
 \def\dj{\hbox{d\kern-.33em\raise.80ex\hbox{-} \raise-.80ex\hbox{\kern-.40em}}}

\begin{document}

\textit{\bf International Journal of Miracles,    01 (2024),pp.
1-39.}

\vspace{3mm}

\begin{center}
\medskip
\epsfig{figure=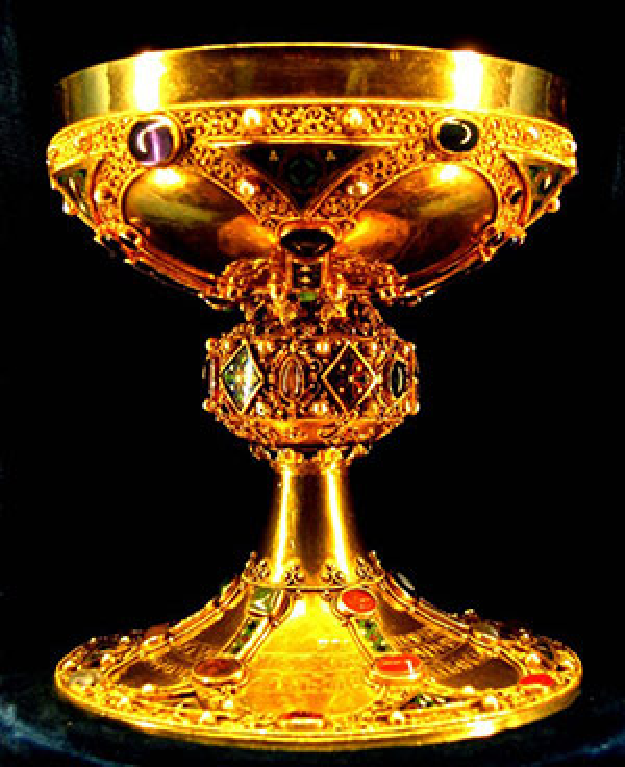,width=5.0cm,height=5.0cm}
\medskip
\footnote{ The 17th century "Holy Grail" of Maths.}
\end{center}
\baselineskip=0.10in

\title[The Lost Proof of Fermat's Last Theorem]{$\mathbb{\quad\;THE\;\;LOST\;\;PROOF\;\;OF\;\;FERMAT'S\;\;LAST\;\;THEOREM}\;\;\;
\mathbb{\;\;\;\;\;\;}$}

\author{$\mathrm{ANDREA\;\;OSSICINI}$}

\maketitle

\begin{center}
{\bf PREMISE}
\end{center}

\vspace{2mm}

\noindent This work contains two papers: the first published in 2022
and entitled "On the nature of some Euler's double equations
equivalent to Fermat's last theorem" provides a marvellous proof
through the so-called discordant forms of appropriate Euler's double
equations, which could have entered in a not very narrow margin and
the second instead published in 2024 and entitled "Some
Diophantus-Fermat double equations equivalent to Frey's elliptic
curve" provides the possible proof, which Fermat has not published
in detail, but which uses the characteristic of all right-angled
triangles with sides equal to whole numbers, or the famous
Pythagorean identity.

\noindent Some explanations in session(III) are provided: the first
makes evident the nature of the "proof a' la Fermat" and the
subsequent sessions clarify the direct and interesting connection of
the two elementary proofs and it is necessary if you want to
understand how two different elementary proofs of Fermat's Last
Theorem are possible.

\noindent Regarding the first paper, a method is used that
simplifies Wiles' theory, a theory that has received much honors
from the entire mathematical community.

\noindent More precisely, through the aid of a Diophantine equation
of second degree solved at first not directly, but as a consequence
of the resolution of the double Euler equations that originated it
and finally in a direct way, the author was able to obtain the
following result: the intersection of the infinite solutions of
Euler's double equations gives rise to an empty set and this only by
exploiting a well-known Legendre Theorem, that is a criterion which
concerns the properties of all the Diophantine equations of the
second degree, homogeneous and ternary.

\noindent It must be observed that this proof  must in no way be
interpreted as a sort of absurd revenge of elementary number theory
over more modern analytic and algebraic treatments.

\noindent The author himself has added a section in which he
connects his concepts with some of those used by Wiles in his
complex demonstration.

\noindent This implies that, to a certain extent, Wiles'
demonstration inspired the author of first paper.

\noindent Ultimately in this paper we will illustrate how only
thanks to some of Euler's discoveries was it possible to shed light
on the so-called too narrow margin never written by Fermat.

\noindent For this reason we will also provide some details on an
article that was the real inspiration for achieving these results
(see Last Conclusions).

\vspace{6mm}

\noindent {I report that Mathematics - An Open Access Journal from
MDPI - published an updated and complete version of the first work,
entitled "On the nature of some Euler's double equations equivalent
to Fermat's last theorem"  (27 Nov. 2022): this document is
available online at the following links:

\vspace{2mm}

\noindent https://www.mdpi.com/2227-7390/10/23/4471/htm

\vspace{4mm}

\noindent https://www.mdpi.com/2227-7390/10/23/4471

\vspace{2mm}

\noindent Furthermore, this version makes it clear that we are not
in the presence of a direct proof of Fermat's Last Theorem, but that
of a reformulation of it which indirectly proves the aforementioned
Theorem.}

\vspace{2mm}

\noindent While the second work, entitled "Some Diophantus-Fermat
double equations equivalent to Frey's curve",has been published in
2024 in "Journal of Ramanujan Society of Mathematics and
Mathematical Sciences"  Vol. 11, No. 1 (2023), pp. 177-187,  by
publisher Ramanujan Society of Mathematics and Mathematical Sciences
(RSMAMS).

\vspace{2mm}

Free access

\vspace{1mm}

\noindent
https//rsmams.org/journals/articleslist.php?volume=11\&issue=1\&tag=jrsmams\&part=0.

\vspace{3mm}

{\bf classification} : 11D41  (primary), 11G05 (secondary).

{\bf keywords} : Fermat's Last Theorem, Arithmetic algebraic
geometry, Diophantine geometry.

\vspace{3mm}

\begin{center}
{\bf DIGRESSION}
\end{center}

\vspace{3mm}

\noindent Fermat wrote that his proof would not fit into the margin
of his copy of Arithmetica, and Wiles's 100 pages of dense
mathematics certainly fulfils this criterion, but surely the
Frenchman did not invent modular forms, the Taniyama-Shimura
conjecture, Galois Groups and the Kolyvagin-Flach method centuries
before anyone else.

\vspace{2mm}

\noindent If Fermat did not have Wiles's proof then what did he
have?

\vspace{2mm}

\noindent Mathematicians are divided into two camps:

\vspace{2mm}

\noindent The sceptics believe that Fermat's Last Theorem was the
result of a rare moment of weakness by the 17th-century genius.

\vspace{2mm}

\noindent They claim that although Fermat wrote, "I have discovered
a truly marvellous proof", he had in fact found only a flawed proof.

\vspace{2mm}

\noindent Other mathematicians, the romantic optimists, believe that
Fermat may have had a genuine proof.

\vspace{2mm}

\noindent Whatever this proof might have been, it would have been
based on 17th-century techniques, and would have involved an
argument so cunning that it has eluded everybody.

\vspace{2mm}

\noindent Indeed, there are plenty of mathematicians who believe
that they can still achieve fame and glory by discovering Fermat's
original proof.

\vspace{2mm}

\noindent In my case it is pure passion for the Mathematics and the
desire to do justice to Fermat and his genius !!!

\begin{center}

{\bf  I)   On the nature of some Euler's double equations equivalent
to Fermat's last theorem}.

\end{center}

\begin{center}
{\bf abstract}
\end{center}

In this work we illustrate that a possible proof of Fermat's Last
Theorem derives from an appropriate use of the concordant forms of
Euler and from an equivalent ternary quadratic homogeneous
Diophantine equation able to accommodate a solution of Fermat's ex-
traordinary equation. Following a similar and almost identical
approach to that of A. Wiles, I tried to translate the link between
Euler's double equations (concordant / discordant forms) and
Fermat's Last Theorem into a possible proof of the Fermat Theorem.
More precisely, through the aid of a Diophantine equation of second
degree, homogeneous and ternary, solved not directly, but as a
consequence of the resolution of the double Euler equations that
originated it, I was able to obtain the following result: the
intersection of the infinite solutions of Euler's double equations
gives rise to an empty set and this only by exploiting a well-known
Legendre Theorem, which concerns the properties of all the
Diophantine equations of the second degree, homogeneous and ternary.
The impossibility of solving the second degree Diophantine equation
thus obtained is certainly possible also through methods known and
discovered by Fermat.

\section{\bf Introduction}

Fermat's last theorem affirms : If $n$ is an integer, greater than
2, there are not any positive integers $X, Y, Z$, so that it can be
valid:
\[
X^n + Y^n = Z^n.
\]
Fermat himself proved it for $n$=4 ([7], pp.
108-112),([6],~II,~Chap. XIII, \S $ $ 202--209); it is consequent
its validity also for $n$ as a multiple of 4, because, if $n$ is
equal $4p$, for some positive integer $p$,
\[ X^n + Y^n = Z^n \quad
\Rightarrow  \quad \left( {X^p} \right)^4 + \left( {Y^p} \right)^4 =
\left( {Z^p} \right)^4 \quad \] and this is impossible.

\vspace{2mm}

\noindent In the same way if we succeed in proving the theorem for a
certain $ k-$exponent, then it is valid for all the multiples of
$k$.

\vspace{1mm}

\noindent As every positive integer greater than 2 is divisible
either by a prime odd number (that is different from 2), or by 4, it
will be then sufficient to prove the theorem for all those cases in
which the exponent is a prime odd number ([9], pp. 203-207).

\vspace{1mm}

 \noindent In this proof we will discuss all those cases
in which the exponent $n$ is an odd $\qquad\qquad$ number $>1$ and,
from now onwards, we will indicate the Fermat Last Theorem with the
acronym F.L.T..

\section{\bf Indeterminate Analysis of Second Degree}

\noindent Our goal is to take care of the resolution, into {\it
integers}, of quadratic equation with {\it integer coefficients},
depending on $n$ unknowns ([1], Cap. I, pp. 60-69).

\noindent We will develop our considerations on the equation in
three unknowns:
\begin{equation}
F\left( {X,Y,Z} \right) = aX^2 + bY^2 + cZ^2 + dXY + eXZ + fYZ =
0\;\end{equation} \noindent warning that, all what we will say,
extends immediately to the case of $n$ unknowns.

\vspace{1mm}

\noindent Since the (1) is an equation homogeneous, if ($A,B,C$) are
the solutions also ({\it mA, mB, mC}) are solutions.

\vspace{1mm}

\noindent Therefore we deem identical two solutions such as
($A,B,C$) and ({\it mA, mB, mC}).

\vspace{1mm}

\noindent Such assumption, will narrow the search to the only
primitive solutions of Eq.(1), that is, to those in which $X, Y$ and
$Z $ are pairwise relatively prime.

\noindent Let ($x, y, z$) be a solution in integers of the Eq.(1)
and then $F\left( {x,y,z} \right) = 0$ and we put:
\begin{equation} X = \rho \cdot x + \xi  , Y = \rho \cdot y +
\eta  , Z = \rho \cdot z + \zeta \end{equation}

\noindent where $\xi ,\eta ,\zeta $ are arbitrary integer constants
and $\rho $ an unknown to be determined, so that Eqs.(2) provide an
integer solution for Eq.(1).

\vspace{2mm}

\noindent It must be: $F\left( {X,Y,Z} \right) = \rho ^2\left[ {ax^2
+ by^2 + cz^2 + dxy + exz + fyz} \right] + $

\[
 + \rho \cdot \left[ {2a\xi \cdot x + 2b\eta \cdot y + 2c\zeta \cdot z +
d\left( {\xi \cdot y + \eta \cdot x} \right) + e\left( {\xi \cdot z
+ \zeta \cdot x} \right) + f\left( {\eta \cdot z + \zeta \cdot y}
\right)} \right] +
\]
\[
 + \;\left[ {a\xi \xi + b\eta \eta + c\zeta \zeta + d\xi \eta + e\xi \zeta +
f\eta \zeta } \right] = 0.
\]

\vspace{2mm}

\noindent But the coefficient of $\rho ^2$, equal to $F\left(
{x,y,z} \right)\;\;$, is null and the known term is $F\left( {\xi
,\eta ,\zeta } \right)\;$, so, set equal to $M$ (with $ M \ne 0$ due
to the arbitrary of $\xi ,\eta ,\zeta $), the coefficient $\rho $ of
the above equation is equal to $\rho = - \frac{F\left( {\xi ,\eta
,\zeta } \right)}{M}\;$.

\vspace{2mm}

\noindent Consequently, if it is known an integer solution of
Eq.(1), we have infinite other, by putting in Eqs.(2), in place of
$\rho $, the value now found; then, without the divisor $M$, we
have:
\begin{equation} X = \xi \cdot M - xF\left( {\xi ,\eta ,\zeta } \right)\; ;
Y = \eta \cdot M - yF\left( {\xi ,\eta ,\zeta } \right) ;
\end{equation}
\[Z = \zeta \cdot M - zF\left( {\xi ,\eta ,\zeta }
\right).\]

\vspace{1mm}

\noindent These are the general solutions of Eq.(1).

\vspace{2mm}

\noindent To prove it, we will show, by appropriately selecting $\xi
,\eta ,\zeta $, the previous solutions provide a solution of Eq.(1),
given arbitrarily.

\vspace{2mm}

\noindent Let this ($A,B,C$),it is meanwhile F($A,B,C$)=0; if now,
in Eqs.(3) we write $\xi =A$,$\eta =B$,$\zeta =C$,we have the
solution: $ X=AM;Y=BM;Z=CM$, that, without the factor $M$, it is
identified with the one already provided.

\vspace{4mm}

\noindent In conclusion:

\vspace{2mm}

\textbf{Theorem 2.1}: \textit{Let (x, y, z) be an integer solution
of Eq.(1)}. \textit{All its integer solutions are given by Eqs.(3),
without the integer divisor $M$}.

\vspace{2mm}

\noindent Now we solve the  equation $F\left( {X,Y,Z} \right) = X^2
+ aY^2 - Z^2 =0 $  in integer numbers.

\vspace{1mm}

\noindent Keeping in mind that this equation is homogeneous we know
that we can consider identical the two solutions, as (1,0,1) and
$\left( {m,0,\,m} \right)$ .

\vspace{1mm}

\noindent Let's consider, at this point, the trivial solution
(1,0,1) and we will have: $M = 2\left( {\xi - \zeta } \right)\;$;
$F\left( {\xi ,\eta ,\zeta } \right) = \xi ^2 + a\eta ^2 - \zeta
^2\;$ for which all the solutions, keeping in mind the Eqs.(3), are
given by the relations:

\vspace{1mm}

\[
X = 2\xi \left( {\xi - \zeta } \right) - \xi ^2 -a\eta ^2 + \zeta ^2
= \;\left( {\xi - \zeta } \right)^2 - a\eta ^2 \quad ; \quad Y =
2\eta \left( {\xi - \zeta } \right)
\]
\[
Z = 2\zeta \left( {\xi - \zeta } \right) - \xi ^2 -a\eta ^2 + \zeta
^2 = - \;\left( {\xi - \zeta } \right)^2 - a\eta ^2.
\]

\vspace{1mm}

\noindent Therefore assumed $\left( {\xi - \zeta } \right) = \theta
$ and observed that from a solution $\left( {x,y,z} \right)$ we get
others changing sign to one, or two, or all $\left( {x,y,z}
\right)$, we have:

\[
X = \theta ^2 - a\eta ^2 \quad ; \quad Y = 2\theta \,\eta \quad ;
\quad Z = \theta ^2 +a\eta ^2
\]

\vspace{1mm}

\noindent which provide us with all the primitive integer solutions
of quadratic equation, without an appropriate integer divisor M.

\noindent In general we have that all integer solutions for the
equation $X^2 + aY^2= Z^2$ are:
\begin{equation}  X = k\left( {\theta ^2 - a\eta ^2} \right);\quad Y = k\left( {2\theta \,\eta
}\right);\quad Z = k\left( {\theta ^2 + a\eta ^2}
\right).\end{equation}

\noindent where $\theta,\eta $ are natural numbers and $k$ a
rational proportionality factor(see also [3], kap. V, \S 29, pp.
39-44).

\section{On Homogeneous Ternary Quadratic Diophantine Equations $aX^2 + bY^2 - cZ^2=0$}

\textbf{Theorem 3.1}:  \textit{Let $x^n + y^n = z^n$, with $(x,y)=1$
and $n \geq 3$ has a solution, then there exists an equation $ax^2 +
by^2 = cz^2$, where $a,b,c$ are relatively prime and reduced to the
minimum terms, whose a solution could be reduced to a solution of
Fermat's equation}.

\vspace{1mm}

{\bf Proof.}

\vspace{1mm}

\noindent Let $ X_1, Y_1, Z_1$  be three whole numbers pairwise
relatively prime such as to satisfy the Fermat equation $x^n +y^n =
z^n$.

\vspace{1mm}

\noindent Then the following homogeneous ternary quadratic
Diophantine equation, with  $(V,T,P)=1$ exists:
\begin{equation}
\label{eq.5} X_1^n V^2  + Y_1^n T^2 = Z_1^n P^2.
\end{equation}

\noindent We observe that with the following particular nontrivial
solutions:

\vspace{1mm}

\noindent $V=1, T=1$  and $ P=1$ or  $V=T=P$ in Eq.(5) we obtain the
fundamental Hypothesis (Reductio ad Absurdum) of the F.L.T.: \[
{X_1}^n + {Y_1}^n = {Z_1}^n .\]

\vspace{1mm}

\noindent Now by the evident solutions, indicated above, we can
derive an infinite number of solutions of Eq.(5).

\vspace{1mm}

\noindent Let's remember that for Legendre's Theorem if a ternary
quadratic  homogeneous Diophantine equation (assuming $a,b$ and $c$
are fixed) has an integral solution, then the number of possible
solutions is infinite.

\newpage

\noindent Having said this, it is possible to transform the previous
Diophantine equation (5) into the following equivalent Diophantine
equation, with $(V',T',P')=1$ :

\vspace{2mm}

\begin{equation}
\label{eq.6} X_1 V'^2  + Y_1 T'^2 = Z_1 P'^2.
\end{equation}

\vspace{1mm}

\noindent It is sufficient to assume $V' = X_1^k V, T'= Y_1^k T,
P'=Z_1^k P$ where $k = \frac{n-1}{2}$ and $n>1$ is an odd number.

\vspace{1mm}

\noindent Using the "fundamental theorem of Arithmetic" we can
represent ([13], Theorem 19, p. 31):

$\qquad\qquad \qquad X_1 = X_0 U_1^2, \; Y_1 = Y_0 U_2^2, \; Z_1 =
Z_0 U_3^2.$

\vspace{1mm}

\noindent In this case is possible to transform the previous
Diophantine equation (6) into the following equivalent Diophantine
equation with the relative coefficients reduced to the minimum
terms:

\[ X_0 V''^2  + Y_0 T''^2 = Z_0 P''^2.\]

\noindent In fact just assume $V'' = U_1 V', T''= U_2 T', P''= U_3
P'$

\noindent We observe that $ X_0, Y_0, Z_0 $  are pairwise relatively
prime and square-free numbers.

\noindent The proof ends here by properly verifying also the nature
of exponent $n$.

\vspace{3mm}

\section{From the Concordant Forms of Euler to Fermat's Last Theorem}

\vspace{2mm}

\noindent Let $m,n \in \emph{\textbf{Z}} $ $\backslash \left\{ 0
\right\}$ be integers with $m \ne n\;$. Following Euler (see
\cite{ref5}), the quadratic forms $X^2 + mY^2 $ and $X^2 + nY^2 $
(or the numbers $m$ and $n$ themselves) are called $concordant$ if
there are integers $ (X,Y,Z,T) $ with $ Y \ne 0 $ such that:
\begin{equation} X^2 + mY^2 = \,Z^2\quad \quad X^2 + nY^2 = \,T^2.
\end{equation}
\noindent In 1780 Euler seeks criteria for the treatment of the
double equations (7) and his interest and our own turns to proofs of
impossibility for the cases $m$=1, $n$=3 or 4 and others equivalent
to these two ([15], Chap. III, \S XVI, pp. 253-254).

\noindent In practice, Euler called $X^2 + mY^2$ and $X^2 + nY^2$
$concordant$ forms if they can both be made squares by choice of
integers $X,\,Y$ each not zero; otherwise, $discordant$ forms.
\noindent At this stage, let us introduce the following Euler double
equations:
\begin{equation}
 \label{eq.8} P^2 + Y_1^n Q^2 = V^2, \; P^2 - X_1^n Q^2 =
\,T^2\end{equation} with $X_1^n + Y_1^n = Z_1^n $ and $n>1$ odd
number.

\noindent By multiplying the first two equations (8) together, and
multiplying by $\frac{P^2}{Q^6}$, with $P \ne 0$ and $Q \ne 0$,we
get([8]):
\begin{equation}
 \label{eq.9} \frac{P^2V^2T^2}{Q^6} = \frac{P^6}{Q^6} + \left( {Y_1^n
- X_1^n } \right)\frac{P^4}{Q^4} - X_1^n Y_1^n
\frac{P^2}{Q^2}.\end{equation}

\noindent If we then replace $\frac{P^2}{Q^2}$  by $\,X$ and also
$\frac{PVT}{Q^3}$  by $\,Y$ we find that

 \[ Y^2 = \,X\left( {X - X_1^n } \right)\left( {X + Y_1^n }
\right).\]

\vspace{2mm}

\noindent This is known as Frey Elliptic curve ([4], pp. 154--156).

\vspace{5mm}

\noindent In Mathematics, a Frey curve or Frey--Hellegouarch curve
is the elliptic curve:
\begin{equation} Y^2 = \,X\left( {X - X_1^n } \right)\left( {X + Y_1^n
} \right)\end{equation} or, equivalently :
\begin{equation} Y^2 = \,X\left[ {X^2 + X\left( {Y_1^n - X_1^n }
\right) - X_1^n Y_1^n } \right]\end{equation}
 \noindent associated with a (hypothetical) solution of Fermat's equation : $ X_1^n + Y_1^n
= Z_1^n.$

\vspace{1mm}

\noindent In fact, the discriminant
\[
\Delta = \sqrt {\left( {Y_1^n - X_1^n } \right)^2 + 4X_1^n Y_1^n \;}
= \,X_{1}^n + Y_1^n = Z_1^n,
\]
that determines the existence of the polynomial
\[ \left( {X - X_1^n } \right)\left( {X + Y_1^n
} \right) = X^2 + X\left( {Y_1^n - X_1^n } \right) - X_1^n Y_1^n\]
\noindent is a perfect power of order n.

\vspace{2mm}

\noindent Frey suggested, in 1985, that the existence of a
non-trivial solution to $\quad$  $X^n + Y^n = Z^n$  would imply the
existence of a non-modular elliptic curve, viz. $Y^2 =
X(X-X^n)(X+Y^n)$. \vspace{1mm}

\vspace{2mm}

\noindent This suggestion was proved by Ribet in 1986.

\vspace{2mm}

 \noindent This curve is semi-stable and in 1993 Wiles
announced a proof (subsequently found to need another key
ingredient, furnished by Wiles and Taylor) that every semi-stable
elliptic curve is modular, the semi-stable case of the
Taniyama-Shimura-Weil conjecture ([16]  and [14]).

\vspace{2mm}

\noindent Hence no non-trivial $X^n + Y^n = Z^n$ can exist.

\vspace{2mm}

\noindent Basically thanks to the spectacular work of A. Wiles,
today we know that Frey's elliptic curve  not exist and from this
derives indirectly, as an absurd, the F.L.T..

\vspace{2mm}

\noindent Now, multiplying the first two equations (8) respectively
by $X_1^n$ and by $Y_1^n$ and at end adding together we get the
following homogeneous ternary quadratic equation (see Section 3):
\begin{equation}
 \label{eq.5}X_1^n V^2 + Y_1^n T^2 = Z_1^n P^2
\end{equation}
with the identity $X_1^n + Y_1^n  = Z_1^n $ and  $n>1$ odd number.

\noindent So, we can also enunciate the following conjecture:

\vspace{4mm}

\textbf{Conjecture 4.1}: \textit{Fermat's Last Theorem is true only
if the homogeneous ternary quadratic Diophantine equation (12) does
not exist (in the sense that the Diophantine equation (12) has no
integer solutions )}.

\vspace{2mm}

\noindent Nobody prevents us from assuming the evident solution
$V=T=P=1$ or $V=T=P$ in the equazion (12) and with this we obtain
the solution of Fermat equation: $X_1^n + Y_1^n = Z_1^n$.

\noindent Now from  the  Euler double equations (8) by subtracting,
we have:
\[ V^2 - T^2 = \,Z_1^n Q^2.\]

\noindent This equation together with equation (12) gives rise to a
system perfectly equivalent to Euler's double equations (8) (see
section 5).

\vspace{3mm}

\noindent We have also with  $V=T=1$ or $V=T$ :

\[ V^2 - T^2 = \,Z_1^n Q^2=0.\]

\noindent By definition, in Euler's $concordant$ forms, $Q$ is
absolutely non-zero integer.

\vspace{1mm}

\noindent It follows that $Z_1^n=0$ and the homogeneous ternary
quadratic Diophantine equation (12) does not exist.

\noindent We observe that the same result can be achieved
immediately if we assume $V=T=P=1$ or $V=T=P$ already in Eqs.(8), in
fact with $Q$ non-zero integer we even have $X_1^n=Y_1^n=0$ and
therefore still $Z_1^n=0$.

\vspace{1mm}

\noindent Further verification of these conclusions is also possible
in this way.

\noindent Let us introduce the following Euler double equations:
\begin{equation}
 P'^2 + Y_1^n Q^2 = V^2, \; P''^2 - X_1^n Q^2 =
\,T''^2\end{equation} with $X_1^n + Y_1^n = Z_1^n $ and  $n>1$ odd
number or
\begin{equation}
 P'^2 + Y_1^n Q^2 = V^2, \; P'''^2 - X_1^n Q'^2 =
\,T'''^2\end{equation} with $X_1^n + Y_1^n = Z_1^n $ and  $n>1$ odd
number.

\noindent From Eqs.(4) we have the following solutions of first
Euler equation of Eqs.(13):
\begin{equation} P'= k\left( {\theta ^2 - Y_1^n\eta ^2} \right),\quad  Q= k\left( {2\theta
\,\eta }\right),\quad  V= k\left( {\theta ^2 + Y_1^n\eta ^2}\right)
\end{equation}
and the following solutions of second Euler equation of Eqs.(13):
\begin{equation} P''= k\left( {\theta ^2 + X_1^n\eta ^2} \right),\quad  Q= k\left( {2\theta
\,\eta }\right),\quad  T''= k\left( {\theta ^2 - X_1^n\eta
^2}\right)
\end{equation}
or the following solutions of second Euler equation of Eqs.(14):
\begin{equation} P'''= k'\left( {\theta '^2 + X_1^n\eta '^2} \right),\; Q'= k'\left(
{2\theta '\,\eta ' }\right),\; T'''= k'\left( {\theta'^2 - X_1^n\eta
'^2}\right).
\end{equation}

\noindent Now assuming  $ V =T =P $  with $Q$ non-zero integer we
have the following result due to Eqs.(15) and Eqs.(16):
\[P = P'= P'' \Rightarrow - Y_1^n= X_1^n \Rightarrow
Z_1^n= 0\ \quad and \]
\[ V = T'' \Rightarrow Y_1^n= - X_1^n \Rightarrow Z_1^n= 0. \]
\noindent While, with Eqs.(15) and Eqs.(17), we have:
\[P = P' = V \Rightarrow - Y_1^n= Y_1^n \Rightarrow
Y_1^n= 0\ \quad and \]
\[P = P''' = T''' \Rightarrow X_1^n= - X_1^n \Rightarrow
X_1^n= 0\] and therefore still $Z_1^n=0$.

\vspace{1mm}

\noindent In conclusion what has been described so far in relation
to Conjecture 4.1 obviously does not have a demonstrative value, but
allows us to state the following equivalent theorem:

\vspace{1mm}

\textbf{Fundamental Theorem}: \textit{Fermat's Last Theorem is true
if and only if is not possible a solution in integers of Eqs.(8)
with $Q$ non-zero integer, that is these are discordant forms.}

\vspace{2mm}

\noindent In practice, this means that if the system of quadratic
Eqs.(8) admits only the trivial solutions $(m,0,\pm m,\pm m)$, that
include also (1,0,1,1), then the quadratic forms $P^2 + Y_1^n Q^2$
and $P^2 - X_1^n Q^2$ are a fortiori called $discordant$.

\vspace{1mm}

\noindent A complete and direct proof of this Theorem is formed in
section 6.

\vspace{1mm}

\section{The Nature of Euler's Double Equations Through the Algebraic Geometry}

\vspace{1mm}

\noindent In this section we will concentrate on the following
Euler's concordant/discordant forms Eqs(8):

\[
 \label{eq.1} P^2 + Y_1^n Q^2 = V^2, \; P^2 - X_1^n Q^2 =
\,T^2\] with $X_1^n + Y_1^n = Z_1^n $ and $n \geq 3$.

\vspace{1mm}

\noindent In determining the nature of the Euler double equations
and of an appropriate equivalent Diophantine system, we will make
use of the description given by A. Weil ([15], Chap. II, App. IV,
pp. 140--149) in order to provide some theoretical background to
Fermat's and Euler's method of descent employed in the treatment of
elliptic curves.

\vspace{1mm}

\noindent For simplicity we consider the case where the roots of a
cubic $\Gamma $ are rational integers $\alpha ,\,\beta $ and $\gamma
$.

\vspace{1mm}

\noindent The cubic $\Gamma $ is then given by
\begin{equation}
 \label{eq.2} y^2 = f\left( x \right) = \left( {x - \alpha }
\right)\left( {x - \beta } \right)\left( {x - \gamma }
\right).\end{equation}

\noindent Weil consider an oblique quartic $\Omega \left( {A,B,C}
\right)$ in the space $\left( {u,v,w} \right)$
\begin{equation}
 \label{eq.3} Au^2 + \alpha = Bv^2 + \beta = Cw^2 + \gamma\end{equation}
with $u,v,w \in \;${\Large Q} and the following mapping of $\Omega $
in $\Gamma $
\begin{equation}
 \label{eq.4} x = Au^2 + \alpha ,\qquad y = \sqrt {ABC}
uvw\end{equation} where $A \cdot B \cdot C$ has to be a {square}.

\vspace{3mm}

\noindent In practice Weil states that the determination of rational
points of the curve $\Gamma $ can be reduced to that of finding
rational points of one or more appropriate quartics, such as (19),
given a set of integers $A,B,C$ (positive or negative), considered
squarefree, that is, not divisible by any square greater than 1, and
such that the product $A \cdot B \cdot C$ is a square.

\vspace{2mm}

\noindent In homogeneous coordinates, $\Omega \left( {A,B,C}\right)$
may be regarded as defined by the equation
\begin{equation}
 \label{eq.5} AU^2 + \alpha T^2 = BV^2 + \beta
T^2 = CW^2 + \gamma T^2,\end{equation}

\noindent with integers $U,V,W,T$ without a common divisor.

\vspace{1mm}

\noindent Subsequently, after affirming that Eq.(21) admits at least
one solution, instead of defining

\noindent $\Omega = \Omega \left( {A,B,C} \right)$ through (19),
Weil writes it through the equation of two quadrics in $P^3$, that
is: $\Phi = \sum\limits_{i,j = 1}^4 {a_{ij} X_i Y_j } $ and  $\Psi =
\sum\limits_{i,j = 1}^4 {b_{ij} X_i Y_j }$, with the condition $\Phi
= \Psi = 0$.

\vspace{8mm}

\noindent In detail, one has:
\[ \Phi \left( {U,V,W,T} \right) = \alpha \left( {\beta - \gamma }
\right)\left( {AU^2 + \alpha T^2} \right) + \beta \left( {\gamma -
\alpha } \right)\left( {BV^2 + \beta T^2} \right) + \]
\[\gamma \left( {\alpha - \beta } \right)\left( {CW^2 + \gamma T^2}
\right) = \alpha \left( {\beta - \gamma } \right)AU^2 + \beta \left(
{\gamma - \alpha } \right)BV^2 + \gamma \left( {\alpha - \beta }
\right)CW^2 - \delta T^2
\]

\vspace{1mm}

$\Psi \left( {U,V,W,T} \right) = \left( {\beta - \gamma }
\right)AU^2 + \left( {\gamma - \alpha } \right)BV^2 + \left( {\alpha
- \beta } \right)CW^2$

\vspace{2mm}

\noindent where one has put $ \;\; \delta = \left( {\beta - \gamma }
\right)\left( {\gamma - \alpha } \right)\left( {\alpha - \beta }
\right)$.

\vspace{2mm}

\noindent With this in mind, we consider the following assumptions
\begin{equation}
 \label{eq.6} A = 1,\quad \alpha = 0,\quad B = 1,\quad \beta = X_1^n ,\quad C = 1,\quad \gamma = - Y_1^n
.\end{equation} \noindent In this case Eq.(18) would be reduced to
the Frey elliptic curve :
\begin{equation}
 \label{eq.7} Y^2 = f\left( X \right) = \left( {X}
\right)\left( {X - X_1^n} \right)\left( {X + Y_1^n}
\right).\end{equation}

\vspace{2mm}

\noindent and the Euler double equations (8) with the following
assumptions, in order: $P = U,\; T=W,\; Q = T$ would be reduced to
the oblique quartic $\Omega \left( {A,B,C} \right)=\Omega \left(
{1,1,1} \right)$:

\begin{equation}
\label{eq.8} U^2 = V^2 - Y_1^n T^2 = W^2 + X_1^n T^2.
\end{equation}

\noindent The product $ABC$ is, as required, a perfect square, and
therefore it is certainly possible the application (19) of the
quartic $\Omega $ on cubic $\Gamma $.

\vspace{2mm}

\noindent The expressions of the two quadrics in $P^3$ become

\[\Phi \left( {U,V,W,T} \right) = -Y_1^n X_1^n V^2 + X_1^n Y_1^n  W^2 + Z_1^n X_1^n Y_1^n T^2\quad\]  and
\[
\Psi \left( {U,V,W,T} \right) = -\left( {Y_1^n + X_1^n } \right) U^2
+ X_1^n V^2 + Y_1^n W^2  = -\left( {Z_1^n} \right) U^2 + X_1^n V^2 +
Y_1^n W.^2
\]

\vspace{2mm}

\noindent Finally, by $\Phi = \Psi = 0$, they are translated into
\vspace{1mm}
\begin{equation}
\label{eq.9} \left( {V^2 - W^2} \right) = \left( {Z_1^n }
\right)T^2\end{equation} \noindent and
\begin{equation}
\label{eq.10} X_1^{n} V^2 + Y_1^{n} W^2= Z_1^n U^2 .
\end{equation}

\vspace{2mm}

\noindent Now Eq.(25) and Eq.(26) with the following replacements:

\[ T \Rightarrow W ,\quad Q \Rightarrow T,\quad P \Rightarrow U \]

\noindent are none other than the equations what we have described
in the section 4, that is:
\[
\label{eq.11} \left( {V^2 - T^2} \right) = \left( {Z_1^n }
\right)Q^2\]
 and
\[
\label{eq.11} X_1^{n} V^2 + Y_1^{n} T^2= Z_1^n P^2 .\]

\noindent This alternative procedure confirms the validity of the
our conclusions: more precisely, I am referring to the fact that
Euler's double equations, as representatives of an evident oblique
quartic of genus 1, can also be defined by means of a pair of
equations of two quadrics in $P^3$, which establish uniquely that
the following Diophantine systems are perfectly equivalent:

\begin{equation}
 \quad \left\{ {\begin{array}{l}
 \; \quad P^2 + Y_1^n Q^2 = V^2\quad \\
  \quad P^2 - X_1^n Q^2 = T^2 \\
 \end{array}} \right.
\quad   \quad \left\{ {\begin{array}{l}
 \;\quad X_1^n V^2 + Y_1^n T^2 = Z_1^n P^2\\
  \quad  Z_1^n Q^2 = V^2 - T^2.\\
  \end{array}} \right.
\end{equation}

\section{The determination of the parameter $Q$ in Euler's double equations}

Let us consider the first Diophantine equation of the second system
(27):
\begin{equation}
X_1^n V^2 + Y_1^n T^2 = Z_1^n P^2
\end{equation}
\noindent and we apply Theorem 2.1.

\noindent Now we solve the equazion $F\left( {X,Y,Z} \right) = aX^2
+ bY^2 - cZ^2 = 0\;$.

\noindent Keeping in mind that this equation is homogeneous we known
that we can consider identical the two solutions, as $(1,1,1)$ and
$(m,m,m)$.

\noindent Let's consider, at this point , the solutions $(1,1,1)$
and we will have:
\[M = 2\left( {a\xi + b\eta - c\zeta } \right)\;; \quad F\left( {\xi
,\eta ,\zeta } \right) = a\xi ^2 + b\eta ^2 - c\zeta ^2\;\]for which
all the solutions, without the integer divisor $M$, keeping in mind
Eq.(3), are given by the relations:
\[
X = a\xi ^2 - b\eta ^2 + 2b\xi \eta + c\zeta \left( {\zeta - 2\xi }
\right) ; \quad Y = - a\xi ^2 + b\eta ^2 + 2a\xi \eta + c\zeta
\left( {\zeta - 2\eta } \right)\]
\[
Z = - a\xi ^2 - b\eta ^2 - \zeta \left[ {c\zeta - 2\left( {a\xi +
b\eta } \right)} \right].
\]

\noindent Without loss of generality, we assume that $\zeta = 0$,
therefore we reduce the intervention of the three integers
 $\xi ,\eta $ and $\zeta $ and to only two of them.

\vspace{5mm}

\noindent In practice we use the following equations instead of
Eqs(2):

\[X = \rho \cdot x + \xi ,\quad Y = \rho \cdot y + \eta ,\quad  Z = \rho
\cdot z\]

\noindent and eliminates the parameter $\rho $ to obtain the
following parametric solutions of Eq.(28):
\begin{equation} V = \lambda \left( {X_1^n \xi ^2 - Y_1^n \eta ^2 + 2Y_1^n \xi \eta }
\right); \; T = \lambda \left( { - X_1^n \xi ^2 + Y_1^n \eta ^2 +
2X_1^n \xi \eta } \right); P = \lambda \left( {X_1^n \xi ^2 + Y_1^n
\eta ^2} \right).
\end{equation}
\noindent Where $\xi $ and $\eta $ are coprime integers and $\lambda
$ is rational proportionality factor.

\vspace{1mm}

\noindent Moreover $\; \xi $ , $\eta $ and $\lambda $ are uniquely
determinated, up to a simultaneous change of sign of $\xi $ and
$\eta $.

\vspace{1mm}

\noindent One standard method of obtaining the above parametrization
can be found also in ([2], \S 6.3.2, pp. 343-346).

\vspace{1mm}

\noindent Now from the second equation of the second system (27)
with the Eqs.(29) and  $\left( {V,T} \right) = 1$, we have with
$\lambda = \frac{1}{M} $ :

\[ Z_1^n Q^2 = V^2 - T^2 = \frac{1}{M^2}\left[ {4\xi \eta \left(
{\xi - \eta } \right)\left( {X_1^n \xi + Y_1^n \eta } \right)\left(
{X_1^n + Y_1^n } \right)} \right] \Rightarrow\]
\begin{equation}
Q^2 = \frac{1}{M^2}\;4\xi \eta \left( {\xi - \eta } \right)\left(
{X_1^n \xi + Y_1^n \eta } \right).
\end{equation}

\noindent For the last factor $\left( {X_1^n \xi + Y_1^n \eta }
\right)$ we can consider the following linear equation:
\begin{equation}\left( {X_1^n \xi + Y_1^n \eta } \right) = hZ_1^n
\end{equation}
\noindent which certainly, admitting the obvious solution $\xi=\eta
=h$, provides us all solutions also with $\xi \ne \eta $, that is:
\begin{equation} \xi = h + Y_1^n \theta;\quad \eta = h - X_1^n \theta.
 \end{equation}
\noindent Besides we have:

\begin{equation}
\left( {\xi - \eta } \right) = Z_1^n \theta .
 \end{equation}

\noindent Therefore bearing in mind that $\left( {X_1 ,Y_1 ,Z_1 }
\right) = 1$, $\left( {V,T,P} \right) = 1$ and $\left( {\xi ,\eta }
\right) = 1$, we have also that $\left( {h,\theta } \right) = 1$.

\noindent Now, Eq.(30) with Eq.(31), Eq.(33) and in addition with $
M = 2\left( {a\xi + b\eta } \right) = 2\left( {X_1^n \xi + Y_1^n
\eta } \right)$ provides:
\begin{equation}
Q^2 = \frac{1}{4\left( {X_1^n \xi + Y_1^n \eta } \right)^2}\;\; 4\xi
\eta \left( {\xi - \eta } \right)\left( {X_1^n \xi + Y_1^n \eta }
\right) = \frac{\xi \eta \theta \,Z_1^n }{h\,Z_1^n } = \xi \eta
\frac{\theta }{h}.
 \end{equation}

\noindent Now we will resort to the Corollary 6.3.8 ([2], p. 346).

\noindent In the case of $\left( {V,T,P} \right)= 1$ we have that
the rational proportionality parameter

\noindent in the Eqs.(29) is $\; \lambda = \frac{1}{r}\;$ with $
r\,|{2Y_1^n Z_1^n} $.

\noindent Now $\lambda = \frac{1}{M}\; \Rightarrow h = \frac{Y_1^n
}{m}\;$ with $m\,\left| {\,Y_1^n } \right.\;$.

\noindent Without loss of generality, we can verify  only the
following extreme case $\;m = 1$ and $\;m = Y_1^n$ [ see Section 8:
Appendix ].

\noindent In fact, thanks to the solutions (32), a single and
appropriate value of  $h$  is sufficient for these equations to
constitute the general solution of the linear equation (31).

\noindent It follows that for $\theta  = 0, \pm 1, \pm  2,... $
formulas (32) give all the integral solutions of equation (31).

\noindent  The necessary condition is that $h$  is an exact divisor
of  $Y_1^n$  and consequently $h=Y_1^n$ or $h=1$ both satisfy this
condition.

\vspace{1mm}

\noindent In the first case with $\;h = Y_1^n $ we have  from
 Eq.(34): $Q^2 = \left( {1 + \theta } \right) \left( \theta \right)
\left( {Y_1^n - X_1^n \theta } \right)$ with the three positive
factors in brackets that are pairwise relatively prime.

\noindent By the uniqueness of the prime decomposition we have
$\left( {1 + \theta } \right) $  and $\,\theta $ should be equal to
squares and this is absurd.

\vspace{2mm}

\noindent In the second case with $h = 1\;,\theta > 0$ and $X_1^n<0$
we have from Eq.(34): $Q^2 = \left( {1 + Y_1^n \theta } \right)
\left( \theta \right) \left( {1 - X_1^n \theta } \right)$ with the
three positive factors in brackets that are pairwise relatively
prime.

\vspace{1mm}

\noindent By the uniqueness of the prime decomposition we have that:
\[ \begin{array}{l}
 \xi = \left(1 + Y_1^n \theta \right)= V_1^2 \,; \eta = \left(1 -
X_1^n \theta \right)= T_1^2 \,;P_1^2=1 \,;\theta = Q_1^2
\,.\end{array}
\]

\noindent In conclusion we have the further double Euler equations:

\[P_1^2 + Y_1^n Q_1^2 = V_1^2 \quad ; \quad P_1^2 - X_1^n Q_1^2 =
T_1^2\]

\noindent with $Q > Q_1$, if compared with the double Euler
equations of the first Diophantine system (27).

\vspace{1mm}

\noindent Repeating the argument indefinitely would the give a
sequence of positive integer $Q > Q_1 > Q_2 > Q_3 > ...$ which
decreased indefinitely.

\vspace{1mm}

\noindent This is impossible, because imply an "infinite descent"
for parameter Q.

\vspace{1mm}

\noindent The determination of the parameter $Q$, as rational
integer not equal to zero, ends here, but we must remember that the
Eq.(34) was determined only thanks by assuming the obvious solution
$\xi=\eta =h$ of the linear equation (31).

\vspace{1mm}

\noindent In this case due to Eq.(33), assuming $Z_1^n > 0$, we have
$\theta = 0 $ and this results in the zeroing of the parameter Q.

\noindent The double equations of Euler are discordant forms and so
the F.L.T. turns out to be true, just as honestly announced by
Fermat himself.

\vspace{2mm}

\section{Conclusions}

\vspace{2mm}

\noindent In this paper we have try to prove F.L.T. making use of
elementary techniques, certainly known to P. Fermat.

\vspace{2mm}

\noindent We show that making use of the concordant forms of Euler
and a ternary quadratic homogeneous Diophantine equation, it is
possible to derive a proof of the F.L.T. without recurring to modern
techniques, but exploiting the important criterion of Legendre for
determining the solutions of ternary quadratic homogeneous equation.

\noindent The proof, here presented, is valid in the case of all odd
exponents greater than one (see the proof of the Theorem 3.1).

\noindent We observe however that also in the  case of exponent $n =
4$ the double equations of Euler are discordant: in this case, in
the double equations of Euler, defined by the expressions (7) is
just assume that $m = -n = 1$ .

\vspace{1mm}

\noindent More precisely we have the following system of equations:
\[\left\{ {\begin{array}{l}
 \;  X^2 + Y^2 = Z^2 \\
 \; X^2 - Y^2 = T^2\\
  \end{array}} \right.
\]
\noindent that has no solutions in the natural numbers.

\vspace{2mm}

\vspace{1mm} \noindent This theorem of a "congruent number" was
anticipated by Fibonacci in his book "The Book of squares" ([12],
Chap. III, \S $ $ VI-2, pp. 310--311), but with a demonstration does
not complete (the first complete proof was provided by Fermat with
the equivalent Theorem: {{\it No Pythagorean triangle has square
area}})([13] ,Chap. II, pp. 50--56).

\vspace{2mm}

\noindent In this work we have not used the proof of non-existence
of the Frey elliptic curve, but we have limited ourselves to proof
of non-existence of the single homogeneous ternary quadratic
equation Eq.(5), defined in the proof of the Theorem 3.1, but whose
origin [see Eq.(12)] is implicit in the nature of Euler's double
equations.

\vspace{1mm}

\noindent The double equations of Euler gave rise in different ways
to the elliptic curve of Frey and to a particular homogeneous
ternary quadratic equation: both characterized by the presence of
${X_1}^n$, ${Y_1}^n$ and ${Z_1}^n$ in their coefficients.

\vspace{1mm}

\noindent For this it was possible to use a similar strategy to
build a proof of the F.L.T..

\vspace{5mm}

{\bf Additional Remarks}

\vspace{3mm}

\indent Remark 1. This work is a reworking of an incomplete essay
$\ll$ Euler's double equations equivalent to Fermat's Last Theorem
$\gg$ ([11]) with the aim of making a proof absolutely complete of
the F.L.T. and consequently making accessible a Theorem of which
Fermat claimed to have a proof and which generations of
mathematicians have tried in vain to try to rediscover it.

\vspace{2mm}

\indent Remark 2. In 1753 Euler calls the Fermat Last Theorem $\ll$
a very beautiful theorem $\gg$, adding that he could only prove it
for $n = 3$ and $n = 4$ and in no other case ([15], Chap. III, \S $
$ 5-d, p. 181).

\vspace{1mm}

\noindent In 1770, He gave a proof with exponent $p$~=~3, in his
{Algebra} ([6],~II,~Chap. XV, \S $ $ 243), but his proof by infinite
descent contained a major gap.

\vspace{1mm}

\noindent However, since Euler himself had proved the lemma
necessary to complete the proof in other work, he is generally
credited with the first proof.

\vspace{1mm}

\noindent The author of this paper has done nothing but complete a
work begun and masterly conducted by Euler himself.

\vspace{1mm}

\noindent For this reason, he considers himself as a co-author of
this proof, but hopes, as shown elsewhere ([10]), that this way of
working can become a normal habit.

\vspace{5mm}

\section{Appendix}

\vspace{4mm}

\noindent Let us consider the following homogeneous linear equation
$ ax + by + cz = 0.$

\vspace{2mm}

\noindent All integer solutions are given by formulas:
\[
x = \frac{k}{\delta }\left( {b\alpha } \right),\;y = \frac{k}{\delta
}\left( {c\beta - a\alpha } \right),\; z = - \frac{k}{\delta }\left(
{b\beta } \right)
\]
\noindent where $k,\,\alpha ,\;\beta \;$  are integers, $\left(
{\alpha ,\;\beta \;} \right) = 1\;$ and $\delta = \left( {b\alpha
,\;c\beta - a\alpha ,\;b\beta } \right)\;.$

\vspace{2mm}

\noindent Having said this, let us consider the equation $X_1^n \xi
+ Y_1^n \eta - Z_1^n h = 0$.

\vspace{2mm}

\noindent We will have the following integer solutions:

\begin{equation} \xi = \frac{k}{\delta }\left( {Y_1^n \alpha }
\right),\;\;\eta = \frac{k}{\delta }\left( {Z_1^n \beta - X_1^n
\alpha } \right),\quad h = \frac{k}{\delta }\left( {Y_1^n \beta }
\right)\end{equation}

\vspace{1mm}

\noindent where $\left( {\alpha ,\;\beta \;} \right) = 1\;$  and
$\delta = \left( {Y_1^n \alpha ,\;Z_1^n \beta - X_1^n \alpha
,\;Y_1^n \beta } \right).$

\vspace{3mm}

\noindent Alongside these we also consider Eqs.(32), that is:

\begin{equation}\xi = h + Y_1^n \theta ,\;\eta = h - X_1^n
\theta .\end{equation}

\noindent Resulting in any case $h\;\left| {\;Y_1^n } \right.$ and
$\left( {\xi ,\;\eta \;} \right) = 1\;$ we have $k = 1$ and $\left(
{h,\;\theta \;} \right) = 1\;$

\noindent Furthermore, in order to determine values for the
parameter $h$, we consider the following equation [see Eq.(34)] :

\vspace{1mm}

\begin{equation}
Q^2 = {\frac{1}{h} \xi \eta \theta }
 \end{equation}

\noindent From Eqs.(35) we have:$\qquad \frac{\xi }{h} =
\frac{\alpha}{\beta} \Rightarrow \beta= 1 $   and

\begin{equation} \frac{\xi }{h} = \alpha. \end{equation}

\vspace{3mm}

\noindent Furthermore, again from Eqs.(35)

\begin{equation} \xi - \eta  = \frac{Y_1^n \alpha }{\delta} - \frac{1}{\delta} \left(Z_1^n - X_1^n \alpha \right) =
\frac{1}{\delta} Z_1^n \left(\alpha - 1 \right) .
\end{equation}

\vspace{5mm}

\noindent From Eqs.(36) we have:

\begin{equation} \xi - \eta  = Z_1^n  \theta. \end{equation}

\noindent The Eq.(39) and Eq.(40)  $\Rightarrow$ \begin{equation}
\theta \delta = \alpha - 1. \end{equation}

\vspace{1mm}

\noindent Now resulting:

\begin{equation} h \delta = Y_1^n \end{equation}

\noindent we also have:$ \qquad  \frac{h}{\theta} =
\frac{Y_1^n}{\alpha - 1} \Rightarrow $

\begin{equation}  h =  \frac{\theta}{\alpha - 1} Y_1^n  \quad  or  \quad
 Y_1^n =  \frac{\alpha - 1}{\theta} h. \end{equation}

\vspace{2mm}

\noindent From Eqs.(36) with Eqs.(43) we obtain

\begin{equation} \xi = h + Y_1^n \theta =  Y_1^n \theta \frac{\alpha}{\alpha - 1}
\qquad ; \qquad  \eta = h - X_1^n \theta =  \theta \left(\frac{Z_1^n
- X_1^n\alpha}{\alpha - 1}\right).
\end{equation}

\noindent From Eq.(37) with Eq.(38) and Eqs.(44) we have:

\[ Q^2 = \alpha \theta^2 \left(\frac{Z_1^n - X_1^n\alpha}{\alpha -
1}\right) = \frac{\theta \alpha} {\alpha - 1} \left( Z_1^n - \alpha
X_1^n \right) \theta .\]

\noindent At end with Eqs.(43) we obtain the following equivalent
equations:

\[ Q^2 = \left( \frac{\xi} {Y_1^n}\right) \; \left( Z_1^n - \alpha X_1^n \right)
\theta \]

\noindent or

\[ Q^2 =  \xi \;\left( \frac{ Z_1^n - \alpha X_1^n} {Y_1^n}
\right) \theta .\]

\vspace{1mm}

\noindent The determination of the parameter $Q$, as rational
integer not equal to zero, ends here.

\vspace{1mm}

\noindent The former equation  $\;\Rightarrow h = Y_1^n \;$ and $\;
\delta = 1$ [ see Eq.(37) and Eq.(42)] and  the latter equation
$\;\Rightarrow \delta = Y_1^n \;$ and $\; h = 1$  [see Eq.(37) and
second formula of Eqs.(35)].

\newpage

\begin{center}
{\bf  II) \bf Some Diophantus-Fermat double equations equivalent to
Frey's elliptic curve}.
\end{center}

\noindent {\bf abstract}: In this work I demonstrate that a possible
origin of the Frey elliptic curve derives from an appropriate use of
the double equations of Diophantus-Fermat and from an isomorphism: a
birational application between the double equations and an elliptic
curve.

\noindent From this origin I deduce a Fundamental Theorem which
allows an exact reformulation of Fermat's Last Theorem.

\noindent A complete proof of this Theorem, consisting of a system
of homogeneous ternary quadratic Diophantine equations, is certainly
possible also through methods known and discovered by Fermat,in
order to solve his extraordinary equation.

\vspace{3mm}

\noindent{\bf 1. The double equations of Diophantus-Fermat and the  Frey elliptic curve} \vspace{.2cm}\\

\setcounter {equation}{0}

\noindent A careful reading of the existing documentation about the
Diophantine problems, reveals that Fermat, and especially Euler,
often used the so-called "double equations" of Diophantus, that is
$\;$ $ax^2+bx+c=z^2;\; a'x^2+b'x+c'=t^2$ with the conditions that
$a$ and $a',$ or $c$ and $c'$ are {squares}.

\noindent These conditions ensure the existence of rational
solutions of the double equations.

\noindent These equations can be written in a more general form as:
\begin{equation}
\label{eq1} ax^2 + 2bxy + cy^2 = z^2 \quad \quad  a'x^2 + 2b'xy +
c'y^2 = t^2.\end{equation}

\noindent Indeed usually both Fermat and Euler considered only the
curves of those forms which have, in the projective space, at least
one "visible" rational point.

\noindent Fermat and Euler derive from few evident solutions an
infinite number of solutions.

\noindent Under this last hypothesis ([5],Chap. II,Appendix III,pp.
135--139) the curve determined by the equations (1) results
isomorphic to the one given by

\begin{equation}
\label{eq2} Y^2 = \,X\;\left[ {\left( {b'X - b} \right)^2 - \left(
{a'X - a} \right)\left( {c'X - c} \right)} \right]\end{equation}

\noindent i.e. an elliptic curve (see also Appendix A).

\noindent In fact, an elliptic curve, which has at least one
rational point, can be written as a cubic $y^2 = f\left( x
\right)\,$, where $f$ is a polynomial of degree 3.

\noindent Given this, we consider the following system, consisting
of a pair of $\ll$double equations$\gg$\begin{equation} \label{eq3}
\left\{ {\begin{array}{l}
 \;\left( 3 \right)_1 \quad \quad X_1^n V^2 + Y_1^n T^2 =
U'^2\quad \quad V^2 - T^2 = W^2 \\
 \;\left( 3 \right)_2 \quad \quad X_1^n W^2 + Z_1^n T^2 =
U'^2\quad \;\,\,W^2 + T^2 = V^2 \\
 \end{array}} \right.\end{equation}

\noindent where $\; X_{1\;} ,Y_{1\;} ,Z_{1\;} $ are integer numbers
(positive or negative), pairwise relatively primes, $\;n> 2$ is a
natural number  and $\;U',V,W,T\;$ are integer variables.

\noindent Applying the isomorphism described by Eq. (2) we obtain,
from the first two equations of the system (3), i.e.  the $\left( 3
\right)_1 $, the elliptic curve
\begin{equation}
\label{eq4} Y^2 = \,X\left( {X - X_1^n } \right)\left( {X + Y_1^n }
\right),\end{equation} and from the other two equations, the $\left(
3 \right)_2 $, the further elliptic curve
\begin{equation}
\label{eq5} Y^2 = \, - X\left( {X - X_1^n } \right)\left( {X - Z_1^n
} \right).\end{equation}

\vspace{5mm}

Combining Eq. (4) and Eq. (5) and using the relation $ X = {X_1^n }
\mathord{\left/ {\vphantom {{X_1^n } 2}} \right.
\kern-\nulldelimiterspace} 2$ one obtains the following identity:
 \begin{equation}\label{eq6} X_1^n + Y_1^n = Z_1^n .\end{equation}

\setcounter {equation}{6}

\noindent Now the elliptic curve (4), together with the identity
(6), is nothing but the Frey elliptic curve ([1], pp.154--156).

\noindent In Mathematics, a Frey curve, or Frey--Hellegouarch curve,
is the elliptic curve:
\begin{equation} Y^2 = \,X\left( {X - X_1^n } \right)\left( {X + Y_1^n
} \right)\end{equation} or, equivalently :\begin{equation} Y^2 =
\,X\left[ {X^2 + X\left( {Y_1^n - X_1^n } \right) - X_1^n Y_1^n }
\right]\end{equation}
 \noindent associated with a (hypothetical) solution of Fermat's equation : $ X_1^n + Y_1^n
= Z_1^n.$

\noindent In fact, the discriminant
\[\Delta = \sqrt {\left( {Y_1^n - X_1^n } \right)^2 + 4X_1^n Y_1^n \;}
= \,X_{1}^n + Y_1^n = Z_1^n,
\]
that determines the existence of the polynomial $ \left( {X - X_1^n
} \right)\left( {X + Y_1^n } \right) = X^2 + X\left( {Y_1^n - X_1^n
} \right) - X_1^n Y_1^n$ \noindent is a perfect power of order n.

\vspace{1mm}

\noindent Frey suggested, in 1985, that the existence of a
non-trivial solution to $X^n + Y^n = Z^n$  would imply the existence
of a non-modular elliptic curve,viz.  $\qquad\qquad$
 \[Y^2 = X(X-X^n)(X+Y^n).\]

\vspace{1mm}

\noindent This suggestion was proved by Ribet in 1986.

\vspace{1mm}

\noindent This curve is semi-stable and in 1993 Wiles announced a
proof (subsequently found to need another key ingredient, furnished
by Wiles and Taylor) that every semi-stable elliptic curve is
modular, the semi-stable case of the Taniyama-Shimura-Weil
conjecture ([6]  and [4]).

\vspace{1mm}

\noindent Hence no non-trivial $X^n + Y^n = Z^n$ can exist.

\vspace{1mm}

\noindent Moreover, as Euler found out, treating similar problems,
regarding algebraic curves of genus 1, the two problems, connected
to curves (4) and (5), are completely equivalent.

\noindent In our case it is simple to verify that the elliptic curve
(5) can be reduced to (4) by the transformation $X \Rightarrow - X +
X_1^n $ and the identity (6).

\vspace{2mm}

{\bf 2. The Diophantine System }\vspace{.2cm}\\

One can reduce the system (3) to the following Diophantine system
\begin{equation}
\label{eq.9} \quad \left\{ {\begin{array}{l}
 \; \quad X_1^n V^2 + Y_1^n T^2 = U'^2\quad \\
  \quad X_1^n W^2 + Z_1^n T^2 =  U'^2 \quad \\
 \quad W^2 + T^2 = V^2. \\
 \end{array}} \right.
\end{equation}

\noindent Our proof of~Fermat's Last Theorem~ consists in the
demonstration that it is not possible a  resolution in whole
numbers, all different from zero, of a system derived from system
(9), but analogous, [see section 4 and system (19)], with integer
coefficients and using integer variables $U,W',T',V'$.

\noindent From the first two equations of the system $\left(9
\right)$ one obtains
\begin{equation}
\label{eq.10} X_1^n V^2 + Y_1^n T^2 = X_1^n W^2 + Z_1^n
T^2.\end{equation}

\noindent Now from Eq.(10) is
\begin{equation}
\label{eq.12}  X_1^n \left( {V^2 - W^2} \right) = \left( {Z_1^n -
Y_1^n } \right)\;T^2. \end{equation}

\noindent  Eq. (11) results in identity (6) if the third equation in
the systems (9), $W^2 + T^2  = V^2$, is satisfied.

\noindent In fact, since this equation is the Pythagorean triangle,
in general, it accepts the following integer solutions, where $p,q$
are natural numbers and $k$ a proportionality factor (the values of
$W$and $T$ are interchangeable if necessary): \[  W = k\left(
{2pq}\right);\quad T = k\left( {p^2 - q^2} \right);\quad V = k\left(
{p^2 + q^2} \right).\]

\noindent We can therefore consider also the primitive integer
solutions with $p,q \in \;${ ${\rm N}$
\begin{equation}
\label{eq.14} W = 2pq\,;\quad \quad T = p^2 - q^2;\quad \quad V =
p^2 + q^2. \end{equation}

\noindent Thus Eq.(11), with $p$ and $q$ relatively prime, of
opposite parity and $p > q >0 $ now is reduced to the identity (6).

\vspace{0.5cm}

{\bf 3.  On Homogeneous Ternary Quadratic Diophantine Equations}
\[ aX^2 + bY^2 - cZ^2=0.\]\vspace{.2cm}\\

\textbf{Theorem 3.1}:  \textit{Let $x^n + y^n = z^n$, with $(x,y)=1$
and $n \geq 3$ has a solution, then there exists an equation $ax^2 +
by^2 = cz^2$, where $a,b,c$ are relatively prime and reduced to the
minimum terms, whose a solution could be reduced to a solution of
Fermat's equation}.

\vspace{2mm}

{\bf Proof.}

\vspace{2mm}

\noindent Let $ X_1, Y_1, Z_1$  be three whole numbers pairwise
relatively prime such as to satisfy the Fermat equation $x^n +y^n =
z^n$.

\vspace{2mm}

\noindent Then the following homogeneous ternary quadratic
Diophantine equation, with  $(V,T,P)=1$ exists:

\begin{equation}
\label{eq.5} X_1^n V^2  + Y_1^n T^2 = Z_1^n P^2.
\end{equation}

\noindent We observe that with the following particular nontrivial
solutions:

\vspace{1mm}

\noindent $V=1, T=1$  and $ P=1$ or  $V=T=P$ in Eq.(13) we obtain
the fundamental Hypothesis (Reductio ad Absurdum) of the F.L.T.: \[
{X_1}^n + {Y_1}^n = {Z_1}^n .\]

\vspace{1mm}

\noindent Now by the evident solutions, indicated above, we can
derive an infinite number of solutions of Eq.(13).

\vspace{1mm}

\noindent Let's remember that for Legendre's Theorem if a ternary
quadratic  homogeneous Diophantine equation (assuming $a,b$ and $c$
are fixed) has an integral solution, then the number of possible
solutions is infinite.

\vspace{1mm}

\noindent Having said this, it is possible to transform the previous
Diophantine equation (13) into the following equivalent
equation,with $(V',T',P')=1$ :

\vspace{1mm}

\begin{equation}
\label{eq.6} X_1 V'^2  + Y_1 T'^2 = Z_1 P'^2.
\end{equation}

\vspace{1mm}

\noindent It is sufficient to assume $V' = X_1^k V, T'= Y_1^k T,
P'=Z_1^k P$ where $k = \frac{n-1}{2}$ and $n>1$ odd number.

\vspace{1mm}

\noindent Using the "fundamental theorem of Arithmetic" we can
represent ([3], Theorem 19, p. 31): $\qquad X_1 = X_2 U_1^2, \; Y_1
= Y_2 U_2^2, \; Z_1 = Z_2 U_3^2.$

\noindent In this case is possible to transform the previous
Diophantine equation (14) into the following equivalent Diophantine
equation with the relative coefficients reduced to the minimum
terms: \[ X_2 V''^2  + Y_2 T''^2 = Z_2 P''^2.\]

\noindent In fact just assume $V'' = U_1 V', T''= U_2 T', P''= U_3
P'.$

\noindent We observe that $ X_2, Y_2, Z_2 $  are pairwise relatively
prime and square-free numbers.

\noindent The proof ends here by properly verifying also the nature
of exponent $n$.

\vspace{1mm}

\noindent {\bf Theorem 3.2} : Let's suppose that $x^n + y^n = z^n$,
with  $n \geq 3$ has a solution, in this case we will have an
equation $ax^2 + by^2 = cz^2$, where $c$ is a square, whose solution
could be reduced to a solution of Fermat's equation.

{\bf Proof.}

\noindent From Theorem 3.1 [see Eq.(13)] we have the following
equation $X_1^n V^2 + Y_1^n T^2 = Z_1^n U^2$ that could be reduced
to a solution of Fermat's equation with $n \geq 3$ odd integer.

\noindent Now multiplying the coefficient $ X_1^n, Y_1^n, Z_1^n$  by
factor $Z_1^n$ we have
\[X_1^n Z_1^n V^2  + Y_1^n Z_1^n T^2 =  Z_1^n Z_1^n U^2\]

\noindent and with  $X_0=X_1 Z_1, \; Y_0=Y_1 Z_1, \; Z_0=Z_1^2$  we
get
\begin{equation} \label{eq.63} X_0^n V^2 + Y_0^n T^2 = Z_0^n U^2
\end{equation}
\noindent and $Z_0^n=\left( {Z_1^n} \right)^2$, that is a square.

\noindent In this case we have also, with  g.c.d.$\left( {X_0,Y_0}
\right) = Z_1$, the equation of Fermat
\begin{equation}\label{eq.64} X_0^n + Y_0^n = Z_0^n.\end{equation}
\noindent The proof ends here.

\vspace{1cm}

{\bf 4. The Lost Proof}\vspace{.2cm}\\

\noindent At this point,  multiplying the Eq.(15) by factoring
quadratic $Z_0^n$ we have
\[ X_0^n Z_0^n V^2 + Y_0^n Z_0^n T^2 = Z_0^n Z_0^n U^2 \]
\noindent and finally with $V' = Z_0^{\frac{n}{2}}V\;,T' =
Z_0^{\frac{n}{2}}T$ we obtain $ X_0^n V'^2 + Y_0^n T'^2 = \left(
{Z_0^n} \right)^2 U^2.$

\noindent That said, let's consider the following double equations
of Diophantus-Fermat, necessary to give rise to the well known
Frey's elliptic curve
\begin{equation}  X_0^n V'^2 + Y_0^n T'^2 =\left( {Z_0^n} \right)^2
U^2=U''^2\;\; ; \;\; V'^2-T'^2=W'^2:
\end{equation}
\[ Y^2 = \,X\left( {X - X_0^n } \right)\left(
{X + Y_0^n } \right).\]

\noindent that together with the identity (16) can be rewritten in
\begin{equation} Z_0^{2n} U^2 = X_0^n V'^2 + Y_0^n T'^2 = X_0^n W'^2 +
Z_0^n T'^2. \end{equation}

\noindent In practice we have rewritten the system (9)  in the
following Diophantine system:
\begin{equation}
\label{eq.9} \quad \left\{ {\begin{array}{l}
 \; \quad X_0^n V'^2 + Y_0^n T'^2 = Z_0^{2n} U^2\quad \\
  \quad X_0^n W'^2 + Z_0^n T'^2 = Z_0^{2n} U^2 \quad \\
 \quad W'^2 + T'^2 = V'^2. \\
 \end{array}} \right.
\end{equation}

\noindent Eqs.(18) give us:
\begin{equation} U^2\left[ {Z_0^n } \right]^2 - V'^2\left[ {Z_0^n } \right] +
W'^2Y_0^n = 0 \end{equation}

or equivalently

\begin{equation} U^2\left[ {Z_0^n } \right]^2 - T'^2\left[ {Z_0^n }
\right] - W'^2X_0^n = 0. \end{equation}

\noindent $Z_0^n $ is a square, so the product of the two roots in
Eq.(20), through the Viete-Girard formulas,  is

$ \qquad\qquad\qquad \; \left[ {Z_0^n } \right]_1 \cdot \left[
{Z_0^n } \right]_2 = \frac{W'^2Y_0^n }{U^2}
 \Rightarrow Y_0^n $, which is a square,

\noindent and in Eq.(21) is

$ \qquad\qquad\qquad \; \left[ {Z_0^n } \right]_1 \cdot \left[
{Z_0^n } \right]_2 = \frac{-W'^2X_0^n }{U^2}
 \Rightarrow -X_0^n $, which is a square.

\noindent These latest results are certainly true only with the
assumption that $ W'^2 $ is non-zero.

\noindent From Theorem 3.1 we have that $ X_1, Y_1,Z_1$ are pairwise
relatively prime and with $Y_0^n = {\Box }\footnote{ The symbol $
{\Box }$ represents an indeterminate square.}= Y_1^n Z_1^n$ and $
X_0^n = - {\Box }= X_1^n Z_1^n$ we obtain: \[Z_1^n = {\Box } \; ; \;
Y_1^n = {\Box }\; ; \; X_1^n = - {\Box }\;.\]

\noindent With this last result, obtained also thanks to the use of
a Pythagorean equation [see Eqs.(17)], one finds also:

\[ \left[ {Z_0^n } \right]_1 \cdot \left[ {Z_0^n } \right]_2 =
\frac{W'^2Y_0^n }{U^2} =  \frac{-W'^2X_0^n }{U^2}.\]

\noindent This gives finally the special solution:
\[ Y_1^n = - X_1^n  \Rightarrow Z_1^n = 0.\]

\noindent Consequently the Diophantine system (19) does not admit
integer solutions.

\noindent A further confirmation of these conclusions comes from
what is reported below.

\noindent Keeping in mind that Eq. (20) and Eq. (21) have arisen
from rewriting the original System (9) into System (19), we have to
consider the various substitutions we have subsequently applied and
in particular by $V' = Z_0^{\frac{n}{2}}V\;,T' = Z_0^{\frac{n}{2}}T$
e $ W' = Z_0^{\frac{n}{2}}W$   ( because of the Pythagorean identity
$ V^2 = T^2 +  W^2 $ ) we can rewrite Eq. (20) and Eq. (21) as
follows:

\begin{equation} \left[ {Z_1^n } \right]^3 \left[ {Z_1^n } \left(U^2 - V^2 \right) +
W^2Y_1^n \right] = 0  \end{equation}

or equivalently

\begin{equation} \left[ {Z_1^n } \right]^3 \left[ {Z_1^n } \left(U^2 - T^2 \right) -
W^2X_1^n \right] = 0. \end{equation}

\noindent At this point, canceling the second factors of the two
products (22) and (23) we have:

$\qquad \left[ {Z_1^n } \left(V^2 - U^2 \right) \right]=  W^2Y_1^n$
$\qquad$  and $\qquad\qquad \left[ {Z_1^n } \left(U^2 - T^2 \right)
\right] = W^2X_1^n. $

\noindent By dividing them among themselves, with  $Z_1^n$ e $W^2$
different from zero, and simplifying we get back the Diophantine
equation original:
\begin{equation}
 X_1^n V^2 + Y_1^n T^2 = Z_1^n U^2.
\end{equation}

\noindent Therefore this equation can exist only on condition that
$Z_1^n$ e $W^2$ are both non-zero.

\noindent In addition, the quantities  $\left(V^2 - U^2 \right)$ and
$\left(U^2 - T^2 \right)$  cannot be null, because the Pythagorean
identity would imply  $W^2=V^2-T^2=0$.

\vspace{1mm}

\noindent Now from Eqs. (20) and (21), considering the sum of the
roots,  through the Viete-Girard formulas, we have:

$ \qquad\qquad\qquad \; \left[ {Z_0^n } \right]_1 + \left[ {Z_0^n }
\right]_2 = \frac{V'^2}{U^2}$

\noindent and

$ \qquad\qquad\qquad \; \left[ {Z_0^n } \right]_1 + \left[ {Z_0^n }
\right]_2 = \frac{T'^2}{U^2}.$

\noindent The two sums must be equal, therefore from $ V'^2 = Z_0^n
V^2 $  and   $T'^2 = Z_0^n T^2$

\noindent we get that if $ Z_0^n$ is different from zero:

\begin{equation} V'^2 = T'^2 \Rightarrow V^2 = T^2 \Rightarrow W^2 = 0.\end{equation}

\noindent In summary, I would like to state that Eqs (20) and (21)
through their roots (see products and relative sums) they only
provide the following result:

\indent $\qquad\qquad\qquad\qquad Z_1^n=0 \quad or \quad W^2=0. $

\noindent which prove, in an absolute way, on the one hand the non
existence of the original Diophantine equation Eq (24) (as it is not
can build) on the other hand the "power" of the Diophantus System
original (9),which includes among its three homogeneous ternary
equations of second degree also a Pythagorean equation (PT).

\noindent The need for such a soluble Pythagorean equation (PT) in
integers, it is fully justified by a proposition stated and proved
by A. Weil, who established the existence of an isomorifism between
some appropriate double equations of Diophantus-Fermat and a certain
elliptic curve, or the existence of a birational application between
the double equations and an elliptic curve.

\vspace{1cm}

{\bf 5. Analytical digressions}\vspace{.2cm}\\

\noindent There is no doubt that the system $\left(19 \right)$,
inspired by system (9), represents a true "lockpick" of the Fermat
Last Theorem.

\noindent Through the former system, keeping in mind always the
possibility of exchanging the role of $X_0 $ and $Y_0 $ into
identity (16), we are able to establish the following Fundamental
Theorem:

{{\it The Fermat Last Theorem is true \underline {{if and only if }}
a solution in integers, all different from zero, of the following
Diophantine system, made of three homogeneous equations of second
degree, with integer coefficients $X_0^n $,$Y_0^n $,$Z_0^n $, where
$n$ is a natural number $> 2$ and with $U,T',V',W'$ integer
indeterminates is not possible.}}

\begin{equation}
\label{eq.9} \quad \left\{ {\begin{array}{l}
 \; \quad X_0^n V'^2 + Y_0^n T'^2 = Z_0^{2n} U^2\quad \\
  \quad X_0^n W'^2 + Z_0^n T'^2 = Z_0^{2n} U^2 \quad \\
 \quad W'^2 + T'^2 = V'^2. \\
 \end{array}} \right.
\end{equation}

\noindent The presence of a Pythagorean equation in this system has
been proved to be essential, not only to connect the most general
Fermat's equation to the supposed Frey's elliptic curve, but to
demonstrate the above indicated Fundamental Theorem (see Section 4)
and at the end to provide also a proof of Fermat's Last Theorem,
using a method of Reductio ad Absurdum.

\vspace{6mm}

{\bf 6. Conclusions}\vspace{.2cm}\\

\noindent In this paper I demonstrate that a possible origin of
Frey's elliptic curve derives from an appropriate use of the
so-called "double equations" of Diophantus-Fermat and from an
isomorphism: a birational application between the double equations
and an elliptic curve.

\noindent This Frey elliptic curve does not exist ([1], pp.
154--156) and from this derives indirectly, as an absurd, the Fermat
Last Theorem.

\vspace{1mm}

\noindent In this work we wanted to emphasize that a proof of the
Fermat Last Theorem can not be separated by the strong links with
the supposed Frey elliptic curve, although this does not mean that
Fermat, in another way, was unable to produce our own proof.

\begin{center}
{\bf Appendix A. Elliptic Curves from Frey to Diophantus}
\end{center}
\noindent In Mathematics, a Frey curve or Frey--Hellegouarch curve
is the elliptic curve:
\begin{equation} \label{eq.58} Y^2 = \,X\left( {X - X_1^n } \right)\left( {X + Y_1^n
} \right)\end{equation} or, equivalently:\begin{equation} Y^2 =
\,X\left[ {X^2 + X\left( {Y_1^n - X_1^n } \right) - X_1^n Y_1^n }
\right]\end{equation} \noindent associated with a (hypothetical)
solution of Fermat's equation : $ X_1^n + Y_1^n = Z_1^n.$

\noindent In the language of Diophantus and of  Fermat, we consider
the following "double equation":
\begin{equation} ax^2 + 2bxy + cy^2 = z^2 \qquad a'x^2 + 2b'xy + c'y^2 = t^2.\end{equation}

\noindent In Weil's Appendix III ([5], Ch. II, pp.135-139) he
established (modulo the existence of a rational point) an
isomorphism between the curve defined by the equations (29) and a
certain elliptic curve defined by:
\[ Y^2=\,X\;\left[ {\left({b'X - b}\right)^2 -
\left({a'X - a}\right)\left({c'X - c}\right)} \right] =\;\]
\begin{equation}
X\;\left[ {\left({b'^2 - a'c'}\right)X^2 + \left( {ca' + ac' - 2bb'}
\right)X - ac + b^2} \right].\end{equation}

\noindent Let's suppose that the first double equation is $ax^2 +
Y_1^n y^2 = z^2$ .

\noindent In this case we have considered the following assumptions
in Eq.(30): $\quad b = 0$ and $c = Y_1^n.$

\noindent Now the coefficient of $X^2$ in Eq.(28) is equal to
coefficient of $X^2$ in Eq.(30): $\left( {b'^2 - a'c'} \right) = 1$
and the coefficient of $X$ and the known term in Eq.(28) are equal
to the ones in Eq.(30):
\begin{equation} \left({ca' + ac' - 2bb'} \right) = Y_1^n - X_1^n ;
\; \; - ac + b^2 = - X_1^n Y_1^n ,\end{equation}

\noindent but with $b = 0$ and $c = Y_1^n $ we have

\vspace{3mm}

$\qquad \qquad \qquad \quad  - ac = - X_1^n Y_1^n \Rightarrow -
aY_1^n = - X_1^n Y_1^n \Rightarrow a = X_1^n $.

\vspace{3mm}

\noindent From the first of Eq.(31) we have $Y_1^n a' + X_1^n c' =
Y_1^n - X_1^n \Rightarrow \;a' = 1,$

 $\qquad \qquad \qquad \quad c'=- 1\; \Rightarrow b' = 0$.

\vspace{1mm}

\noindent With these results we have the following double equation
of Diophantus:

\noindent $X_1^n x^2+Y_1^n y^2=z^2$ and $x^2-y^2=t^2$ equivalently
to $X_1^n V^2+ Y_1^n T^2=U'^2$ and $\; V^2 - T^2 = W^2\quad $ [see
the equations of the system (3), i.e.  the $\left(3 \right)_1$].

\vspace{3mm}

{\bf Additional} {\bf Remarks}

\vspace{3mm}

\indent REMARK 1.  Fermat's idea, in my opinion, to prove his Last
Theorem, could take place through the following logical steps:

\underline {{\bf 1}}{\bf - } Define a quadratic and homogeneous
ternary equation, in the normal form of Lagrange, able to
accommodate a solution, with n greater than or equal to 3, of its
extraordinary equation (6) [see Theorem 3.2].

\underline {{\bf 2}}{\bf - } Connect this appropriate Diophantine
equation of 2nd degree to the classic Pythagorean equation [see
Eqs.(17)] to build a complete Diophantine system capable of
determining its possible whole solution [see system (19)] .

\underline {{\bf 3}}{\bf - } Establish that this Diophantine system
does not admit congruent integer solutions and therefore as a
consequence of this, there are no three integers that satisfy
Fermat's equation (6).

\vspace{3mm}

\indent REMARK 2. The truth is that the impossibility to solve
single equations can be proved as deduction from the impossibility
of solving a system of equations.

\vspace{1mm}

\noindent The Fundamental Theorem is a reformulation of the Fermat
Last Theorem: his following statements are equivalent:

\vspace{1mm}

\noindent {(A) Fermat's Last Theorem is true }$\; \Leftrightarrow \;
${(A') The Diophantine System does not allow integer solutions
different from zero.}

\noindent Let $n > 2$; there is a bijection between the following
sets:

\vspace{1mm}

\noindent(S) the set of solutions $(x,y,z)$ of Fermat's Equation,
where $x,y,z$  are nonzero natural numbers; and

\vspace{1mm}

\noindent(S') the set of solutions $(u,t',v',w')$ of the Diophantine
System, where $u,t',v',w'$ are nonzero natural numbers.

\vspace{1mm}

\noindent The set of solutions of (S) and (S') are the same, that
gives rise to an empty set, as shown in the Fundamental Theorem.

\vspace{1mm}

\noindent In the literature there are other Diophantine equations,
that were compared to Fermat's equation, i.e. a first result, due to
Lebesgue in 1840, is the following Theorem:

\vspace{1mm}

{\it If Fermat's Last Theorem is true for the exponent $  n  \ge 3$
then the equation $X^{2n} + Y^{2n} = Z^2$  has only trivial
solutions}.

\vspace{1mm}

\noindent The proof of this theorem is extremely simple and is found
in [2].

\vspace{1mm}

\noindent In this case, however, it cannot be said that Lebesgue's
theorem is equivalent to Fermat's Last Theorem, while on the
contrary, the Fundamental Theorem is just equivalent to Fermat's
last theorem.

\vspace{3mm}

\indent REMARK 3. I conclude this work with the following
observation by A. Weil ([5], Chap. IV, \S $ $ VI, pp. 335--336): "
Infinite descent a' la Fermat depends ordinarily upon no more than
the following simple observation: if the product $\alpha \cdot
\beta$ of two ordinary integers (resp. two integers in an algebraic
number-field) is equal to an m-th power, and if the g.c.d. of
$\alpha$ and $\beta$  can take its values only in a given finite set
of integers (resp. of ideals), then both $\alpha$ and $\beta$ are
m-th powers, up to factors which can take their values only in some
assignable finite set." (See  the section 4: The Lost Proof.)

\vspace{1cm}

\begin{center}
{\bf References}
\end{center}
\begin{enumerate}
\item [{[1]}] H. Davenport, {\it The Higher Arithmetic - an introduction to the
theory of number}, 8th ed., Cambridge University Press, New
York,2008.
\item [{[2]}] V.A. Lebesgue, {{\it  Sur un th\`{e}or\'{e}me de Fermat}}{, J. Math. Pures Appl. }{{\bf 5 }}{(1840), pp. 184--185.}
\item [{[3]}] W. Sierpinski, {\it Elementary Theory of Numbers, }Elsevier Science
Publish'ers B.V., Amsterdam,  vol. {\bf 31,}{\it }$\;2^a$ English
ed. 1988.
\item [{[4]}] R. Taylor, A. Wiles, {\it Ring-theoretic properties of certain Hecke
algebras}, Ann. of Math. {\bf 141} (1995), pp. 553--57
\item [{[5]}] A. Weil, {\it Number Theory: an Approach Through History from
Hammurapi to Legendre,} reprint of 1984 Edition , Birkh\"{a}user,
Boston, 2007{\small .}
\item [{[6]}] A. Wiles, {\it Modular elliptic curves and Fermat's Last Theorem},
Ann. of Math. {\bf 141} (1995), pp. 443--551.
\end{enumerate}

\newpage

\vspace{16mm}

\begin{center}
{\bf  III) \bf Pierre De Fermat's secret margin revealed by Leonhard
Euler.}
\end{center}

\maketitle

\vspace{8mm}

\begin{center}
{\bf PREMISE}
\end{center}

\vspace{4mm}

\noindent Recently two elementary proofs of Fermat Last Theorem has
been given by Andrea Ossicini.

\vspace{2mm}

\noindent Both articles effectively provide a reformulation of
Fermat's Last Theorem (F.L.T.).

\vspace{2mm}

\noindent The first, entitled   "On the Nature of Some Euler's
Double Equations Equivalent to Fermat's Last Theorem", has been
published in 2022 in the journal "Mathematics"  by publisher MDPI
(Multidisciplinary Digital Publishing Institute).

\vspace{2mm}

\noindent The Journal "Mathematics" is indexed in SCOPUS. Impact
factor 2.4. It is quoted with a journal rank: JCR - Q1 (Mathematics)
/ CiteScore 3.5 - Q1 (General Mathematics).

\vspace{2mm}

\noindent Ossicini's article is indicated by Mathematics as "Feature
Paper".

\vspace{2mm}

\noindent This label is used to represent the most advanced
investigations which can have a significant impact in the field.

\vspace{2mm}

\noindent A Feature Paper should be an original contribution that
involves several techniques or approaches, provides an outlook for
future research directions and describes possible research
applications. Feature papers are submitted upon individual
invitation or recommendation by the scientific editors and must
receive positive feedback from the reviewers.

\vspace{2mm}

\noindent The second, entitled "Some Diophantus-Fermat double
equations equivalent to Frey's curve",has been published in 2024 in
"Journal of Ramanujan Society of Mathematics and Mathematical
Sciences"  by publisher Ramanujan Society of Mathematics and
Mathematical Sciences (RSMAMS).

\vspace{2mm}

\noindent This Journal is  indexed ZbMath, MathSciNet, EBSCO, ICI,
etc.

\vspace{4mm}

\normalsize

{\bf 1. FERMAT'S LAST THEOREM $X_0^n+Y_0^n=Z_0^n$ and
PYTHAGOREAN EQUATION $V'^2 = W'^2 + T'^2$ . }\vspace{.2cm}\\

\noindent In the work "Some Diophantus-Fermat double equations
equivalent to Frey's curve" we are able to establish the following
{\bf Fundamental Theorem} [3]:

\vspace{3mm}

{{\it The Fermat Last Theorem is true \underline {{if and only if }}
a solution in integers, all different from zero, of the following
Diophantine system, made of three homogeneous equations of second
degree, with integer coefficients $X_0^n $,$Y_0^n $,$Z_0^n $, where
$n$ is a natural number $> 2$ and with $U,T',V',W'$ integer
indeterminates is not possible.}}

\vspace{2mm}

\[
\label{eq.9} \quad \left\{ {\begin{array}{l}
 \; \quad X_0^n V'^2 + Y_0^n T'^2 = Z_0^{2n} U^2\quad \\
  \quad X_0^n W'^2 + Z_0^n T'^2 = Z_0^{2n} U^2 \quad \\
 \quad W'^2 + T'^2 = V'^2. \\
 \end{array}} \right.
\]

\vspace{2mm}

\noindent The presence of a Pythagorean equation in this system has
been proved to be essential, not only to connect the most general
Fermat's equation to the supposed Frey's elliptic curve, but to
demonstrate the above indicated Fundamental Theorem (see Section 4
[3]) and at the end to provide also a proof of Fermat's Last
Theorem, using a method of Reductio ad Absurdum.

\vspace{3mm}

\noindent From the previous Diophantine System we obtained the
following equivalent biquadratic equations (20) and (21) ( see
Section 4 [3]) or see  Eq(1) and Eq(6) in this Section.

\vspace{3mm}

\noindent Bearing in mind the validity of the Pythagorean
equation~$V'^2 = W'^2 + T'^2$ ,  the g.c.d.$\left( {W',T',V'}
\right) = Z_0^{\frac{n}{2}}$  and $X_0^n+Y_0^n=Z_0^n$ with $X_0=X_1
Z_1, \; Y_0=Y_1 Z_1, \; Z_0=Z_1^2$  , we start by the first
biquadratic equation (20):

\begin{equation}
\label{eq.1} U^2\left[ {Z_0^n } \right]^2 - V'^2\left[ {Z_0^n }
\right] + W'^2Y_0^n = 0.\end{equation}

\vspace{3mm}

\noindent Consequently we have that the roots of this equation of
2\r{ }, with unknown~$Z_0^n $, are
\begin{equation}
\label{eq.2} \left[ {Z_0^n } \right]_{1,2} = \frac{V'^2\pm \sqrt
{V'^4 - 4U^2W'^2Y_0^n } }{2U^2}.\end{equation}

\vspace{3mm}

\noindent The expression under the square root must however be a
perfect square, that we put~$N^2$.

\vspace{3mm}

\noindent We have therefore
\begin{equation}
 \label{eq.3} N^2 + 4U^2W'^2Y_0^n = V'^4. \end{equation}

\vspace{3mm}

\noindent Owing to the fact that the two solutions in (2) must be
equal to a square and that their product [(see Eq. (1)] is equal to:
$ \; \left[ {Z_0^n } \right]_1 \cdot \left[ {Z_0^n } \right]_2 =
\frac{W'^2Y_0^n }{U^2}
 \Rightarrow Y_0^n $, which is a square, the Eq.(3) can be represented by the
following Diophantine equation:

\vspace{3mm}

\begin{equation}
 \label{eq.27} x^2 + y^2 = z^4.\end{equation}

\vspace{3mm}

\noindent Thanks to Euler and its Algebra ([1],~II,~Chap. XII, \S \S
$ $ 198-199~) it is possible to obtain the following integer
solutions of Eq.(4):

\vspace{3mm}

$ \quad x = \pm k\left( p^4 - 6p^2q^2 + q^4 \right)= \pm k\left[
\left( {p^2 - q^2} \right)^2 - \left( {2pq} \right)^2\\\right];\;y =
k(4p^3q - 4pq^3) =k\left[ 4pq\left( {p^2 - q^2} \right)\right];\;z
=k( p^2 + q^2).$

\newpage

\noindent Now, taking into account the validity of the Pythagorean
equation $V'^2 = W'^2 + T'^2$, Eq.(3) provides  $N = k\left[\left(
{p^2 - q^2} \right)^2 - \left( {2pq} \right)^2\right]$ and then with
$W' = k(2pq)\;;\;T' = k(p^2 - q^2)\;;\;V' = k(p^2 + q^2)$, Eq.(2)
gives us :

\vspace{3mm}

\begin{equation}
 \label{eq.5} \left[ {Z_0^n } \right]_{1,2} = \frac{V'^2\pm
\left( {T'^2 - W'^2} \right)}{2U^2}.\end{equation}

\vspace{3mm}

\noindent Since Eq.(2) does provide rational solutions  we get:

\vspace{3mm}

\[ \left[{Z_0^n }
\right]_1 = T'^2/U^2.\]

 and

\[\left[{Z_0^n } \right]_2 = W'^2/U^2\]

\vspace{3mm}

\noindent Now through the Viete-Girard formulas we have:

\vspace{3mm}

$\qquad\qquad\qquad  \; \left[ {Z_0^n } \right]_1 \cdot \left[
{Z_0^n } \right]_2 = \frac{W'^2Y_0^n }{U^2} \quad \Rightarrow \left[
{Z_0^n } \right]_1 = Y_0^n \quad \Rightarrow Z_1^n = Y_1^n \quad
\Rightarrow X_1^n = 0.$

\vspace{6mm}

\noindent Bearing again in mind the validity of the Pythagorean
equation~$V'^2 = W'^2 + T'^2$ ,  the g.c.d.$\left(
{W',T',V'}\right)= Z_0^{\frac{n}{2}}$ and $X_0^n+Y_0^n=Z_0^n$ with
$X_0=X_1 Z_1, \; Y_0=Y_1 Z_1, \; Z_0=Z_1^2$  , the Eq.(1) can be in
the equivalent form, that  is :
\begin{equation}
 \label{eq.6} U^2\left[ {Z_0^n } \right]^2 - T'^2\left[ {Z_0^n }
\right] - W'^2X_0^n = 0.\end{equation}

\vspace{2mm}

\noindent Consequently we have that the roots of this equation of
2\r{ }, with unknown~$Z_0^n $, are
\begin{equation}
 \label{eq.7} \left[ {Z_0^n } \right]_{1,2} = \frac{T'^2\pm
\sqrt {T'^4 + 4U^2W'^2X_0^n } }{2U^2}.\end{equation}

\vspace{2mm}

\noindent The expression under the square root must however be a
perfect square, that we put~$N^2$.

\vspace{1mm}

\noindent We have therefore

\begin{equation}
 \label{eq.8} N^2 - 4U^2W'^2X_0^n = T'^4.\end{equation}

\vspace{1mm}

\noindent Owing to the fact that the two solutions in (7) must be
equal to a square and that their product [(see Eq. (6)] is equal to:
$ \; \left[ {Z_0^n } \right]_1 \cdot \left[ {Z_0^n } \right]_2 =
\frac{-W'^2X_0^n }{U^2}
 \Rightarrow -X_0^n $, which is a square, the Eq.(8) can be represented by the
following Diophantine equation

\vspace{1mm}

\begin{equation}
 \label{eq.27} x^2 - y^2 = z^4.\end{equation}

\vspace{1mm}

\noindent Thanks to Euler and its Algebra ([1],~II,~Chap. XII, \S \S
$ $ 198-199~) it is possible to obtain the following integer
solutions of Eq.(9):

\vspace{1mm}

$\quad  x = \pm k\left( p^4 + 6p^2q^2 + q^4 \right)= \pm k \left[
\left( {p^2 + q^2} \right)^2 + \left( {2pq} \right)^2\\\right];\;y =
k(4p^3q + 4pq^3) = k\left[4pq\left( {p^2 + q^2} \right)\right];\; z
= k(p^2 - q^2).$

\vspace{1mm}

\noindent From here, taking into account the validity of the
Pythagorean equation $V'^2=W'^2+T'^2$, Eq.(8) provides:

\vspace{1mm}

$\quad N = k\left[\left( {p^2 + q^2} \right)^2 + \left( {2pq}
\right)^2\right]$ and with $W' =k(2pq)\;;\;T' =k(p^2 - q^2)\;;\;V'=
k(p^2 + q^2)$,

\vspace{1mm}

\noindent Eq.(7) gives us:
\begin{equation}
\label{eq.10} \left[ {Z_0^n } \right]_{1,2} = \frac{T'^2\pm \left(
{V'^2 + W'^2} \right)}{2U^2}.\end{equation}

\vspace{1mm}

\noindent Since Eq.(10) does provide rational solutions, we get

\vspace{1mm}

\[\left[{Z_0^n } \right]_1 = V'^2/U^2\]
  and
\[\left[{Z_0^n } \right]_2 = -W'^2/U^2\]

\vspace{1mm}

\noindent Now through the Viete-Girard formulas we have:

\vspace{1mm}

$\qquad\qquad\qquad  \; \left[ {Z_0^n } \right]_1 \cdot \left[
{Z_0^n } \right]_2 = \frac{-W'^2(X_0^n) }{U^2}=0 \quad \Rightarrow
\left[ {Z_0^n } \right]_1 = X_0^n \quad \Rightarrow Z_1^n = X_1^n
\quad \Rightarrow Y_1^n = 0.$

\vspace{1mm}

\noindent From the roots of the two biquadratic equations, the first
from the equation (1) and the second from equation (6) :

\[ \left[{Z_0^n }
\right]_1 = T'^2/U^2 =  \left[{Z_0^n } \right]_1 = V'^2/U^2\]

\noindent while on the other equality:

\[ \left[{Z_0^n }
\right]_2 = W'^2/U^2 =  \left[{Z_0^n } \right]_2 = -W'^2/U^2.\]

\vspace{1mm}

\noindent At end we obtain : $W'^2=V'^2 - T'^2=0.$

\vspace{1mm}

\noindent All these results are confirmed also by the Viete-Girard
formulas relating to the sum of the roots of equations (1) and (6).

\vspace{1mm}

\noindent In fact we have:

\vspace{1mm}

$ \qquad\qquad\qquad \; \left[ {Z_0^n } \right]_1 + \left[ {Z_0^n }
\right]_2 = \frac{V'^2}{U^2}$

\vspace{1mm}

\noindent and

\vspace{1mm}

$ \qquad\qquad\qquad \; \left[ {Z_0^n } \right]_1 + \left[ {Z_0^n }
\right]_2 = \frac{T'^2}{U^2}.$

\newpage

\noindent The two sums must be equal, therefore from $ V'^2 = Z_0^n
V^2 $  and   $T'^2 = Z_0^n T^2$, we get that if $ Z_0^n$ is
different from zero:

\begin{equation} V'^2 = T'^2 \Rightarrow V^2 = T^2 \Rightarrow W^2 = 0.\end{equation}

\vspace{1mm}

\noindent The proof of Fermat's Last Theorem is complete in that the
diophantine system of the Fundamental Theorem does not admit integer
solutions all other than zero for its indeterminates and its link
with the Pythagorean relation is extraordinary and indispensable for
an "elementary" proof, thanks also to the solution found in {\bf
Euler's Algebra}.

\vspace{5mm}

{\bf 2. Fermat's Last Theorem Proof (Fermat and Euler). }\vspace{.2cm}\\

\setcounter {equation}{0}

\noindent We start from the equation (15) of work:

\vspace{1mm}

"Some Diophantus-Fermat double equations equivalent to Frey's
Elliptic Curve" [3]:

\begin{equation} X_0^n V^2 + Y_0^n T^2 = Z_0^n U^2\end{equation}

\vspace{1mm}

\indent with $\qquad\qquad\qquad\;$ $Z_0^n=(Z_1^n)^2$,  $X_0^n +
Y_0^n = Z_0^n $ e $V^2-T^2=W^2.$

\vspace{1mm}

\indent  Now with the equation (17) of the work [3] and the
application of A. Weil's proposition to construct an isomorphism
with the Frey's elliptic curve we have, also exploiting the normal
form of Lagrange:

\begin{equation}  X_0^n V'^2 + Y_0^n T'^2 = (Z_0^n)^2 U^2=(Z_0^n)U'^2=U''^2\end{equation}

\vspace{2mm}

\indent with $\qquad\qquad\;$  $V'^2-T'^2=W'^2$ ; $V'^2= (Z_0^n)V^2$
; $T'^2= (Z_0^n)T^2$ ; $ U'^2= (Z_0^n)U^2 $.

\vspace{5mm}

\indent  That said, considering the paper  "On the Nature of Some
Euler's Double Equations Equivalent to Fermat's Last Theorem" [2],
where the following theorem was proved:

\vspace{3mm}

[{\bf Fundamental Theorem}] {Fermat's Last Theorem is true if and
only if a solution is not possible in integers of Euler's Double
Equations $P^2 + Y_1^n Q^2 = V^2$ and $P^2 - X_1^n Q^2 = T^2.$  with
the $Q$ non-zero integer; that is, these are discordant~forms.}

\vspace{2mm}

Therefore we can can write the following doubles Euler equations due
to equation (2) of this Section:

\vspace{2mm}

$\qquad\qquad\qquad\qquad $ $U'^2 + Y_0^n Q^2 = V'^2$ $\qquad$  and
$\qquad$ $U'^2 - X_0^n Q^2 = T'^2.$

\vspace{3mm}

\indent From these we immediately derive the following equation:

\vspace{2mm}

$\qquad\quad$ $ Z_0^n Q^2 = V'^2-T'^2= (Z_0^n) (V^2-T^2)=(Z_0^n)
W^2$ which implies $ Q^2=W^2.$

\newpage

\noindent However, keeping in mind that equation (2) of this Section
gives rise to equations (20) and (21) of work [3] due to identity
Pythagorean, that is:

\begin{equation} U^2\left[ {Z_0^n } \right]^2 - V'^2\left[ {Z_0^n } \right] +
W'^2Y_0^n = 0 \end{equation}

or equivalently

\begin{equation} U^2\left[ {Z_0^n } \right]^2 - T'^2\left[ {Z_0^n }
\right] - W'^2X_0^n = 0. \end{equation}

\vspace{1mm}

\noindent applying the Viete-Girard formulas relating to the sum of
the roots of a polynomial we have: $ \qquad \left[ {Z_0^n }
\right]_1 + \left[ {Z_0^n } \right]_2 = \frac{V'^2}{U^2}$ e $ \left[
{Z_0^n } \right]_1 + \left[ {Z_0^n } \right]_2 = \frac{T'^2}{U^2}.$

\vspace{2mm}

\noindent Ultimately we get: $V'^2 = T'^2 $ that is $W^2=V^2-T^2=0$,
which results in $Q = 0$  and therefore doubles Euler's equations
give rise to discordant forms.

\vspace{1mm}

\noindent In this way we have produced a further elementary
verification of Fermat's Last Theorem putting the two together
proofs, that is the proof attributable to Fermat [3] and that
attributable to Euler [2], without using the parametric solutions of
the equation Diophantine (1) of this Section.

\vspace{5mm}

{\bf 3. Digressions on the elementary proofs of the F.L.T.}\vspace{.2cm}\\

\setcounter {equation}{0}

\noindent  Starting from "Some Diophantus-Fermat double equations
equivalent to Frey's Elliptic Curve" [3] and from the proof of
Fermat's Last Theorem, which is based on the non-existence of an
appropriate Diophantine equation, ternary and homogeneous of second
degree, capable of accommodating a possible integer solution of the
Fermat Equation, we recall that we have deduced the following
equations with $Z_0^n=\left( {Z_1^n} \right)^2$ (see Eq.(20) e
Eq.(21) in [3]):

\vspace{1mm}

\begin{equation} U^2\left[ {Z_0^n } \right]^2 - V'^2\left[ {Z_0^n } \right] +
W'^2Y_0^n = 0 \end{equation}

\vspace{1mm}

or equivalently

\begin{equation} U^2\left[ {Z_0^n } \right]^2 - T'^2\left[ {Z_0^n }
\right] - W'^2X_0^n = 0. \end{equation}

\vspace{1mm}

\noindent The Viete-Girard formulas relating to the product and the
sum of the roots of the biquadratic equation gives a extraordinary
result.

\vspace{1mm}

\noindent In fact, if the sum of the roots implies that $ W'^2=0$ ,
then the product of the same, with

\vspace{1mm}

\noindent $W'^2= Z_0^n W^2=0$ implies that the indeterminate $ W^2$
is zero if one of the solutions ( $ \left[{Z_0^n } \right]_1 or
\left[{Z_0^n } \right]_2 $)  is different from zero.

\vspace{1mm}

\noindent To fully understand the consequences of this result it is
necessary to first consider the expressions (1) and (2) like
polynomials.

\vspace{1mm}

\noindent In this case, it is fundamental the role of {\bf the
identity principle of polynomials}, which is based on the
normalization of individual expressions (same coefficients relating
to variables of the same degree).

\vspace{1mm}

\noindent In fact, comparing the various coefficients of the
parametric equations (1) and (2) in this Section we obtain directly:

\vspace{1mm}

$\qquad\qquad\qquad\qquad\qquad\qquad\qquad\qquad$  $ U^2 = U^2 $
per  $\left[ {Z_0^n } \right]^2$

\vspace{1mm}

$\qquad\qquad\qquad\qquad\qquad\qquad\qquad\qquad$ $ V'^2 = T'^2$
per $\left[ {Z_0^n } \right]$

\noindent which implies \[ W'^2  = 0 \]

\vspace{1mm}

\noindent and equating the known terms we obtain:

\[  W'^2Y_0^n  = - W'^2X_0^n\] or again:

\[Z_0^n W'^2=(Z_0^n)^2 W^2 = 0 \]

\noindent which implies \[ W'^2  = 0 \]

\noindent or even in the case of $ W^2  > 0 \Rightarrow Z_0^n = 0
\Rightarrow Z_1^n = 0.$

\noindent If the cancellation of the indeterminate $ W'^2 $ leads
with it the cancellation of $Z_1^n$ we will have that F.L.T. is in
fact established.

\noindent If we resort to Euler's double equations [2] we have the
following equivalent Diophantine systems:

\begin{equation}
 \quad \left\{ {\begin{array}{l}
 \; \quad U'^2 + Y_0^n Q^2 = V'^2\quad \\
  \quad U'^2 - X_0^n Q^2 = T'^2 \\
 \end{array}} \right.
\quad   \quad \left\{ {\begin{array}{l}
 \;\quad X_0^n V'^2 + Y_0^n T'^2 = Z_0^n U'^2\\
  \quad  Z_0^n Q^2 = V'^2 - T'^2.\\
  \end{array}} \right.
\end{equation}

\noindent  From these we deduce that if $Z_0^n = 0$ we have that
$Q^2= \frac{W'^2}{Z_0^n}= \frac{W^2 Z_0^n}{Z_0^n} $ would result
indeterminate, just as happened in the case where Euler's double
equations are actually only one [which can result only for $m=n$,
see Eq.(7) in [2], to be read as $ Y_0^n = - X_0^n $ which confirm
$Z_0^n = 0$], otherwise, i.e. with $Z_0^n > 0$ [assumed hypothesis
for which it is believed that the F.L.T. is false], we will even
have $Q^2= W^2=0$, i.e. the Euler forms are definitively discordant
and consequently F.L.T. it is proven.

\vspace{3mm}

{\bf 4. The direct link between the two elementary Proofs. }\\

\noindent From Eq(1) and Eq(2) of  the previous session, with $ U'^2
= Z_0^n U^2$  we have:

\begin{equation} U'^2\left[ {Z_0^n } \right] - V'^2\left[ {Z_0^n } \right] +
W'^2Y_0^n = 0 \end{equation}

\noindent or equivalently

\begin{equation} U'^2\left[ {Z_0^n } \right] - T'^2\left[ {Z_0^n }
\right] - W'^2X_0^n = 0. \end{equation}

\vspace{2mm}

Equations (4) and (5) together constitute the following oblique
quartic $\Omega \left( {A,B,C} \right)=\Omega \left(
{Z_0^n,Z_0^n,Z_0^n} \right)$ and in homogeneous coordinates, $\Omega
\left( {A,B,C}\right)$ may be regarded as defined by the equation

\vspace{2mm}

\[ AU'^2  = BV'^2 + \beta W'^2 = CT'^2 + \gamma
W'^2,\]

\vspace{2mm}

 \noindent with integers $U',V',T',W'$ and $\beta=Y_0^n$
and $ \gamma=X_0^n$:
\[  Z_0^n U'^2 = Z_0^n V'^2 - Y_0^n W'^2 = Z_0^n T'^2 + X_0^n W'^2. \]

\vspace{2mm}

\noindent Considering that even the double Euler equations can be
represented by an evident oblique quartic of genus 1, that is
$\Omega \left( {A,B,C} \right)=\Omega \left( {1,1,1} \right)$, we
can pose the following condition :

\vspace{2mm}

\[ V'^2-T'^2 = W'^2  =  Z_0^n Q^2\]

\vspace{2mm}

\noindent and the Eq(4) and Eq(5) provide:

\begin{equation}
\left[ {Z_0^n } \right]  \left[ {U'^2 - V'^2 + Q^2  Y_0^n}
\right]=0\end{equation}

\begin{equation}
\left[ {Z_0^n } \right]  \left[ {U'^2 - T'^2 - Q^2  X_0^n}
\right]=0\end{equation}

\noindent Such products are null a condition of having that the
systems (3) are valid, with  ${Z_0^n }>0$.

\vspace{2mm}

\noindent It is clear that we are in the presence of double
equations of Euler [note that Systems (3) are precisely the solved
Systems in work [2] in section 6].

\vspace{3mm}

\noindent The direct link between the two elementary Proofs of the
Fermat Last Theorem is obvious and has been correctly established
\vspace{4mm}

{\bf 5. Additional remarks. }\\

\noindent In this section, some additional logical and mathematical
remarks on significant aspects of my proof in "On the nature of some
Euler's double equations equivalent to Fermat's last theorem" will
be presented.

\vspace{2mm}

\noindent These considerations are closely related to the
demonstration in the form in which it was published.

\vspace{2mm}

\noindent However, they also add some mathematical details which are
not indispensable for the completeness of the proof, but which
reveal more deeply its logical structure.

\vspace{2mm}

\noindent Starting from the system of Euler's double equations, it
is possible to obtain a homogeneous quadratic Diophantine equation
in four indeterminates, which connects Fermat's equation \(X^{n} +
Y^{n} = Z^{n}\ \)and the aforementioned Euler double equations.

\vspace{2mm}

\noindent For, let us consider the first equation of the System of
Euler's double equations written in the two equivalent forms:

\setcounter {equation}{0}

\begin{equation}
 \quad \left\{ {\begin{array}{l}
 \; \quad P^2 + Y_1^n Q^2 = V^2\quad \\
  \quad P^2 - X_1^n Q^2 = T^2 \\
 \end{array}} \right.
\quad   \quad \left\{ {\begin{array}{l}
 \;\quad X_1^n V^2 + Y_1^n T^2 = Z_1^n P^2\\
  \quad  Z_1^n Q^2 = V^2 - T^2.\\
  \end{array}} \right.
\end{equation}

\vspace{3mm}

\noindent This system is equivalent to the second one in (1), namely

\vspace{2mm}

 \(\left\{
\begin{matrix}
X_{1}^{n}V^{2} + Y_{1}^{n}T^{2} = Z_{1}^{n}P^{2} \\
V^{2} - T^{2} = Z_{1}^{n}Q^{2} \\
\end{matrix} \right.\ \), with the condition that Fermat's equation be
solvable.

\vspace{2mm}

\noindent Let us take into the first equation of this system.

\vspace{2mm}

\noindent Subtract and add the quantity \(X_{1}^{n}Y_{1}^{n}Q^{2}\)
(which is different from 0), from the first member ({\bf Ossicini's
Trick).}

\vspace{2mm}

\noindent The following equation in the unknowns \emph{V, T, Q, P}
is thus obtained:

\vspace{1mm}

\begin{equation}
X_{1}^{n}V^{2} + Y_{1}^{n}T^{2} - X_{1}^{n}Y_{1}^{n}Q^{2} +
X_{1}^{n}Y_{1}^{n}Q^{2} = Z_{1}^{n}P^{2}
\end{equation}

\vspace{1mm}

that is
\begin{equation}
X_{1}^{n}{(V}^{2} - Y_{1}^{n}Q^{2}) + Y_{1}^{n}(T^{2}{+
X}_{1}^{n}Q^{2}) = Z_{1}^{n}P^{2}.\end{equation}

\noindent Through the equations
\begin{equation}
P^{2} = V^{2} - Y_{1}^{n}Q^{2},\quad P^{2} = T^{2}{+ X_{1}^{n}Q}^{2}
\end{equation}

in (1), Equation (3) is transformed into Fermat's equation $
X_{1}^{n} + Y_{1}^{n} = Z_{1}^{n}.$

\vspace{2mm}

\noindent As a consequence of this result, if Euler double equations
are discordant, then Fermat's equation cannot admit integer
solutions different from zero.

\vspace{2mm}

\noindent For it is not possible to verify the evident solution of
equation (3) (whose form is (\emph{x,x,x}), reducible to (1,1,1))
for the values indicated in (4), from which an infinite number of
other solutions could be derived.

\vspace{2mm}

\noindent Therefore, the following link connects Ossicini's
Fundamental Theorem to his conjecture:

\vspace{3mm}

\textbf{Fundamental Theorem}: Fermat's Theorem is true if and only
if a solution in integers of Euler equations of the system (1) is
not possible, being \emph{Q} and integer different form 0. This
means that the equations express a discordant form.

\vspace{2mm}

\noindent This theorem, if true, therefore also implies the validity
of the following proposition:

\vspace{3mm}

\textbf{Conjecture}: Fermat's Theorem is true only if the first
homogeneous quadratic ternary Diophantus system \(\left\{
\begin{matrix}
X_{1}^{n}V^{2} + Y_{1}^{n}T^{2} = Z_{1}^{n}P^{2} \\
V^{2} - T^{2} = Z_{1}^{n}Q^{2} \\
\end{matrix} \right.\ \) does not have integer solutions.

\vspace{2mm}

\noindent Therefore, the quaternary quadratic equation (2) allows us
to reach a direct link between the conjecture and the fundamental
theorem. This link was already conspicuous in the published proof,
but the considerations here expounded make the link between
conjecture and theorem even clearer. As a matter of fact, Equation
(3) represents a remarkable integration between Fermat's equation
and Euler's double equations.

\vspace{4mm}

{\bf 6. Last Conclusions. }\\

\noindent Back in 1952, Prof. Umberto Bini (of the school of
Francesco Severi), in relation to Fermat's Last Theorem [0] , stated
(see below Fig.1, in Italian) :

\vspace{2mm}

" As is well known FERMAT reading the Commentaria in Diophantum by
C. G. BACHET DE MEZIRIAC, had made a habit of annotating them in the
margin.

\vspace{2mm}

Concerning the eighth Diophantum problem, which requires and
provides the resolution in rational numbers of the equation
$X^2+Y^2=Z^2$, FERMAT postulate : <<On the contrary, it is
impossible to divide a cube into the sum of two cubes, a fourth
power into two fourth powers, and, in general, any power of degree
greater than two, into two powers of the same degree.

\vspace{2mm}

 I discovered an admirable demonstration of this general theorem that
this margin is too small to contain ".

\vspace{2mm}

\noindent Reflecting on the earnest and honest manner in which the
magistrate, at the Parliament of Toulouse presented to the
mathematical world his arithmetical discoveries, all of which,
sooner or later, turned out to be true, there is no reason to
believe that the Theorem in question is the result of his rash
assertion or bluff.

\vspace{2mm}

\noindent On only one occasion, to my knowledge, did he say
something later found to be incorrect about the series  : $
{{2}^{2}}^{h} +1 $ whose terms reputed to be all prime, while such
were the first five.

\vspace{2mm}

\noindent And then if it is to be believed that FERMAT has really
found a general dimostration of the impossibility in non-zero
integers primes of $X^n + Y^n + Z^n = 0 $, it is not clear why those
who are dealing with the question today consider, in order to
overcome its difficulties, to have recourse to new concepts, to new
theories in the fields of modern higher analysis.

\vspace{2mm}

\noindent It seems to me that one wants to resort to the hydrogen
bomb, whereas perhaps {\bf a skillful shot with  an arquebus} is
enough; one of those shots of which {\bf the great Euler was a
master}.

\vspace{2mm}

\noindent For this reason I recommend careful reading of the
following works:

\vspace{2mm}

1) "On the Nature of Some Euler's Double Equations Equivalent to
Fermat's Last Theorem" [2].

\vspace{2mm}

2) "Some Diophantus-Fermat double equations equivalent to Frey's
curve" [3].

\vspace{5mm}

\medskip
\epsfig{figure=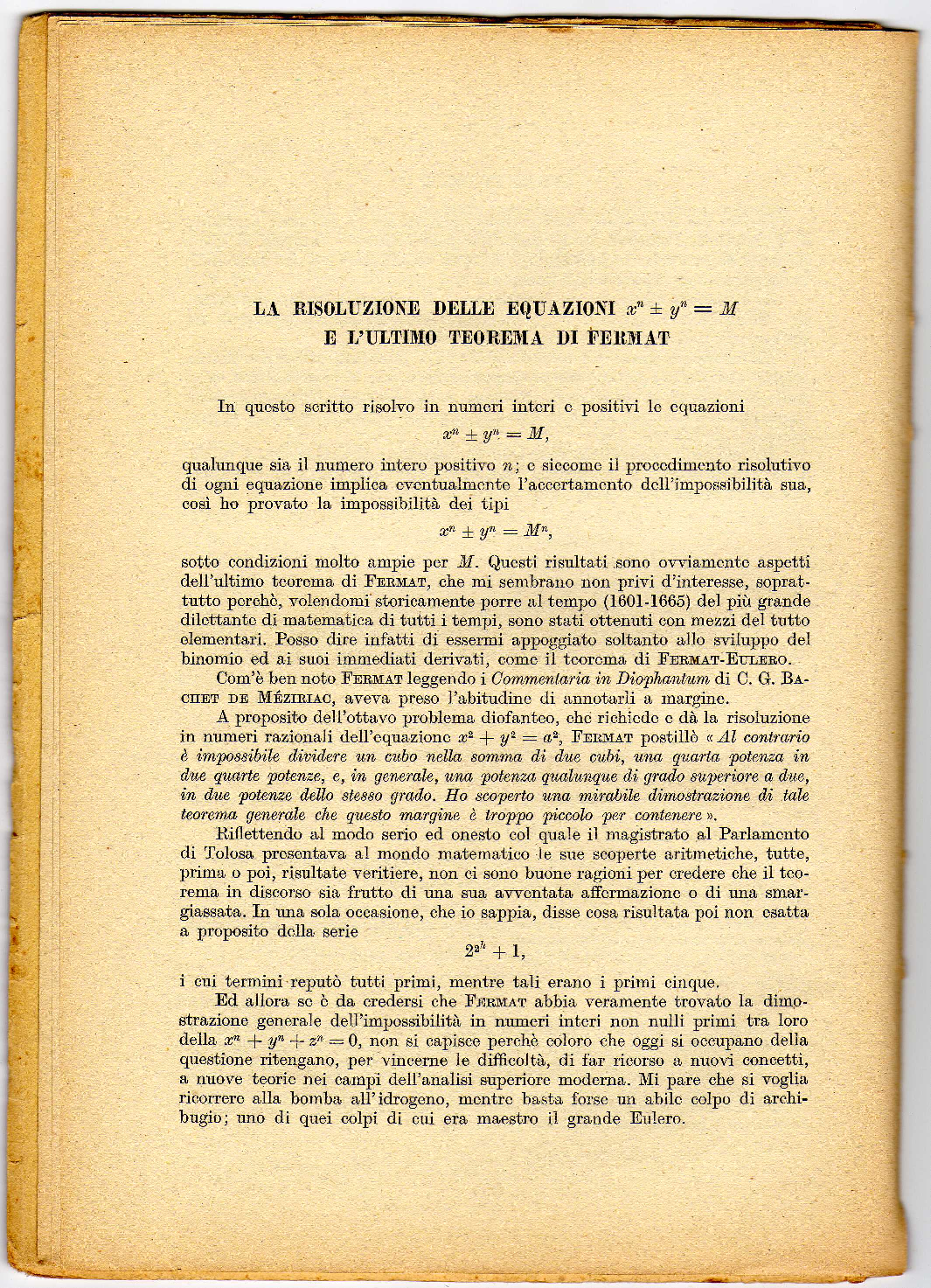,width=14.5cm,height=17.5cm}
\medskip

\vspace{5mm}

{\bf Fig.1  La Risoluzione dell'equazioni $X^n \pm Y^n=M$  e
L'ultimo Teorema di Fermat, ARCHIMEDE}

\newpage

{\bf Dedication.} These papers are written honor of \textbf{Leonhard
Euler}, one of the greatest mathematicians and mechanicians of all
time.

\medskip
\centerline{\epsfig{figure=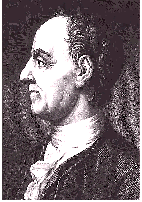,width=12.0cm,height=16.5cm}}

\vspace{2mm}

{\qquad \qquad \qquad \qquad \quad \bf LEONHARD EULER (1707-1783)}

\vspace{2mm}

{\qquad \qquad \qquad \qquad \qquad (from the portrait by A.
Lorgna)}

\medskip
\epsfig{figure=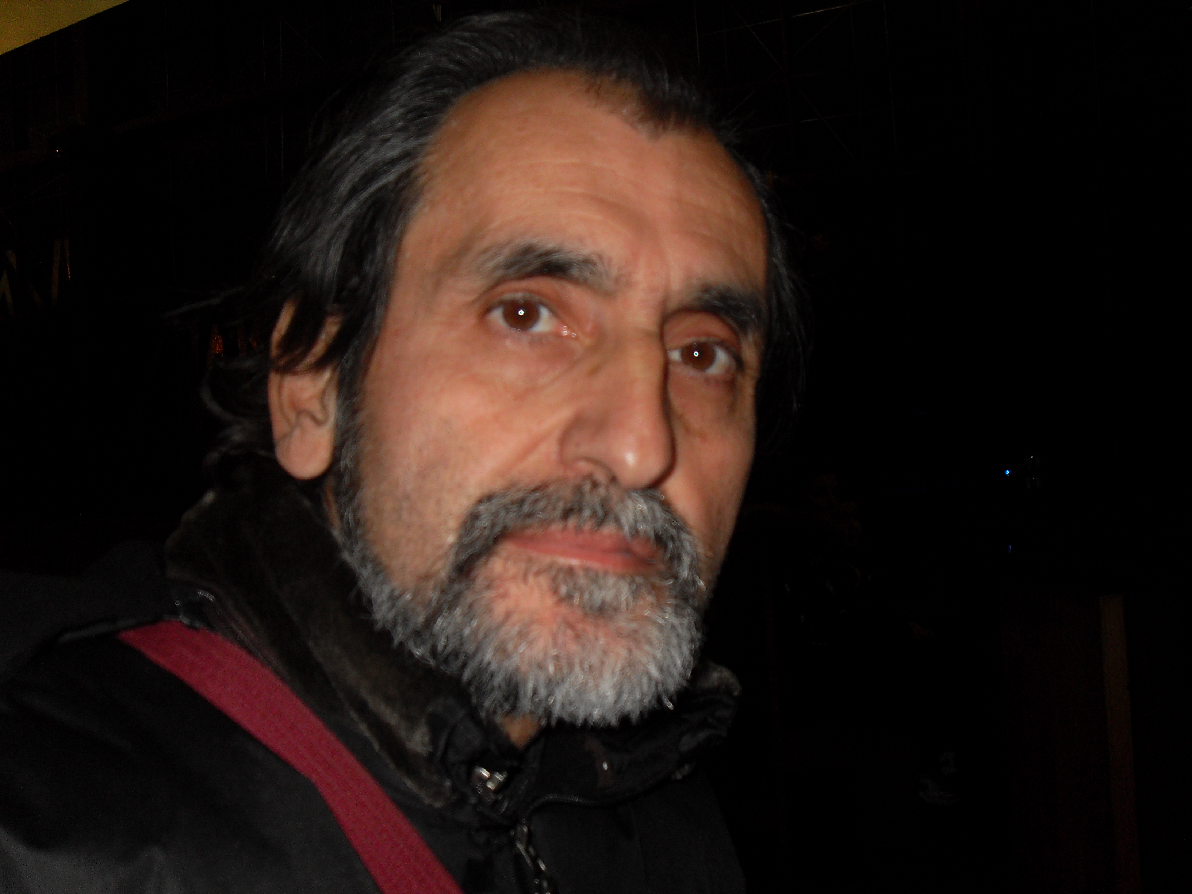,width=8.5cm,height=6.0cm}
\medskip
\baselineskip=0.10in

\vspace{3mm}

\textit{\bf ANDREA OSSICINI}.

\vspace{3mm}

{\bf email} :andrea.ossicini@yahoo.it or andrea.ossicini@gmail.com

\vspace{3mm}

{\bf address} :Via delle Azzorre 352-D2, 00121 Roma, Italy

\end{document}